%
%

\magnification=1200

\font\titfont=cmr10 at 12 pt

\font\headfont=cmr10 at 12 pt



\overfullrule=0in

\def\boxit#1{\hbox{\vrule
 \vtop{%
  \vbox{\hrule\kern 2pt %
     \hbox{\kern 2pt #1\kern 2pt}}%
   \kern 2pt \hrule }%
  \vrule}}

  \def\harr#1#2{\ \smash{\mathop{\hbox to .3in{\rightarrowfill}}\limits^{\scriptstyle#1}_{\scriptstyle#2}}\ }

 \def\Ga{G\aa rding}

\def\J{J}  
\def\bbf{{\bf F}}

\def\jt{\J^2}
\def\jtx{\jt_x}

\def\Edge{{\rm Edge}}

\def\ss{\subset}

\def\half{\hbox{${1\over 2}$}}
\def\smfrac#1#2{\hbox{${#1\over #2}$}}
\def\oa#1{\overrightarrow #1}

\def\log{{\rm log}}
\def\Hess{{\rm Hess}}

\def\tr{{\rm tr}}
\def\max{{\rm max}}

\def\span{{\rm span\,}}

\def\det{{\rm det}}

\def\Sym{{\rm Sym}^2}

\def\rn{\bbr^n}

\def\Int{{\rm Int}}

\def\Symn{{\Sym(\rn)}}

\def\Theorem#1{\medskip\noindent {\bf THEOREM \bf #1.}}
\def\Prop#1{\medskip\noindent {\bf Proposition #1.}}
\def\Cor#1{\medskip\noindent {\bf Corollary #1.}}
\def\Lemma#1{\medskip\noindent {\bf Lemma #1.}}
\def\Remark#1{\medskip\noindent {\bf Remark #1.}}
\def\Note#1{\medskip\noindent {\bf Note #1.}}
\def\Def#1{\medskip\noindent {\bf Definition #1.}}

\def\Ex#1{\medskip\noindent {\bf Example \bf    #1.}}

\def\pf{\medskip\noindent {\bf Proof.}\ }
\def\qed{\hfill  $\vrule width5pt height5pt depth0pt$}
\def\mathqed{\vrule width5pt height5pt depth0pt}
\def\equdef{\buildrel {\rm def} \over  =}

   \def\cp{{\cal P}}

\def\cp{{\cal P}}
\def\cf{{\cal F}}

\def\vf{\varphi}

\def\wt{\widetilde}

\def\and{\qquad {\rm and} \qquad}

\def\bbr{{\bf R}}\def\bbh{{\bf H}}
\def\bbc{{\bf C}}

\def\bbf{{\bf F}}

\def\bbm{{\bf M}}

\def\a{\alpha}
\def\b{\beta}
\def\d{\delta}
\def\e{\epsilon}

\def\g{\gamma}

\def\l{\lambda}
\def\o{\omega}

\def\s{\sigma}

\def\L{\Lambda}
\def\G{\Gamma}
\def\O{\Omega}

\def\bo{\partial \Omega}

\def\Symn{\Sym(\rn)}
 
\def\USC{{\rm USC}}
\def\fa{{\rm\ \  for\ all\ }}

\def\cpt{\wt{\cp}}
\def\ft{\wt F}
\def\ob{\overline{\O}}

\def\Fa{{\oa F}}

\def\Or{{\bbr_+^m}}
\def\hy{hyperbolic }
\def\DG{Dirichlet-G\aa rding }
\def\DGP{\DG polynomial}
\def\ua{{\uparrow}}
\def\lu{\l^{\ua}}
\def\wu{w^{\ua}}
\def\E{E}
\def\Edge{{\rm Edge}}

\centerline{\titfont HYPERBOLIC POLYNOMIALS AND }
\medskip

\centerline{\titfont  THE  DIRICHLET PROBLEM}
\bigskip

\centerline{\titfont F. Reese Harvey and H. Blaine Lawson, Jr.$^*$}
\vglue .9cm
\smallbreak\footnote{}{ $ {} \sp{ *}{\rm Partially}$  supported by
the N.S.F. }

\vskip .5in
\centerline{\bf ABSTRACT} \bigskip
  \font\abstractfont=cmr10 at 10 pt
{{\parindent= .53in
\narrower\abstractfont \noindent

This paper presents a simple, self-contained account of \Ga's theory of 
hyperbolic polynomials, including a recent convexity result  of Bauschke-Guler-Lewis-Sendov
and an inequality of Gurvits.   This account also  contains 
 new results, such as the existence of a real analytic arrangement of
the eigenvalue functions.

In    a second, independent part of the paper, 
the relationship of \Ga's theory  to the authors' recent work on the Dirichlet problem for fully nonlinear partial differential equations  is investigated.
Each \Ga\ polynomial $p$  of degree  $m$ on $\Symn$ (hyperbolic with respect to the identity)
 has an associated 
eigenvalue map $\l:\Symn\to \bbr^m$, defined modulo the permutation group acting on 
$\bbr^m$.  Consequently, each closed symmetric set $E\ss \bbr^m$ induces a second-order 
p.d.e. by requiring, for a $C^2$-function $u$ in $n$-variables, that 
$$
\l\left ((D^2u)(x)\right ) \ \in\ \partial E \fa x.
$$
  Assume  that $A\geq0 \Rightarrow \l(A)\geq0$ and that
$E+\Or\ss E$.  A main result is that for smooth domains $\Omega \ss\rn$ whose
boundary is suitably $(p,E)$-pseudo-convex, the Dirichlet problem has  a unique continuous solution
for all continuous boundary data. This applies in particular to each of the $m$ distinct
branches of the equation $p\left (D^2u\right ) =0$

In the authors' recent extension of results from euclidean domains
to domains in riemannian manifolds, a new  global ingredient, 
called a monotonicity subequation, was introduced.  It is shown in this paper
that for every  \  polynomial $p$ as above, the associated \Ga\ cone
is a monotonicity cone for all branches of the the equation
$p(\Hess\, u) = 0$ where $\Hess \, u$ denotes the riemannian Hessian of $u$.

}}

\vfill\eject

\centerline{\bf TABLE OF CONTENTS} \bigskip

{{\parindent= .1in\narrower\abstractfont \noindent

\qquad 1. Introduction.\medskip

\qquad 2.       G\aa rding's Theory of Hyperbolic Polynomials. \smallskip
\medskip

\qquad 3.  $\G$-Monotone Sets.
\medskip

\qquad 4.    The Dirichlet Problem.\smallskip
 
\medskip

\qquad 5.  Subequations Determined by Hyperbolic Polynomials.

\medskip

\qquad 6. The Dirichlet Problem on Riemannian Manifolds.

\medskip

\qquad Appendix A.  An Algebraic Description of the Branches.

\medskip

\qquad Appendix B.  Gurvits' Inequality.

\medskip

\vfill\eject


\centerline{\headfont \ 1.\  Introduction}
\medskip

This paper is concerned with the G\aa rding theory of hyperbolic polynomials 
and its relation to the Dirichlet problem for fully nonlinear partial differentials equations.

We recall that a homogeneous polynomial of degree $m$ on a finite dimensional
real vector space $V$ is called {\sl hyperbolic with respect to a  direction} 
$a\in V$ if  (we assume  $p(a)>0$)  the one-variable polynomial $t\mapsto p(ta+x)$
has $m$ real roots for each $x\in V$.  Thus we can write
$$
p(ta+x)\ =\ p(a)\prod_{k=1}^m (t+\l_k(x)).
\eqno{(1.1)}
$$
The functions $\l_k(x)$, called the {\sl $a$-eigenvalues of } $x$, are well defined
up to permutation. In 1959 G\aa rding developed a beautiful theory of these
hyperbolic polynomials [G$_2$].  An important part of the theory concerned 
the set 
$$
\G \ \equiv\  \{x\in V : \l_k(x)>0, \ {\rm for \  all\ } k\}
\eqno{(1.2)}
$$
which is proved to be a convex cone.  Of crucial importance for the paper is the {\sl
monotonicity property}:  that if one orders the eigenvalues $\lu_1(x)\leq\cdots\leq \lu_m(x)$, then
$$
\lu_k(x+b) \ > \ \lu_k(x) \fa b\in \G, x\in V.
\eqno{(1.3)}
$$
G\aa rding's  motivation  came from the theory of 
hyperbolic partial differential equations [G$_1$].  (See also [H].) However, in subsequent years his 
 results have attracted attention and found application in  other areas.
 For example, \Ga\  theory was used to initiate hyperbolic programming by 
 Guler [Gu].  (See also [R].)

This paper has two distinct and essentially independent objectives.

The first is to present a simple, self-contained account of G\aa rding's theory.   Our  account starts with a 
new result, Theorem 2.9 below, which states that
\medskip

{\sl 
 For   $b\in \G$ and $x\in V$ the $a$-eigenvalues can be arranged so that each of the functions
 $t\mapsto \l_k(x+tb)$ is a strictly  increasing real-analytic mapping from $\bbr$ onto $\bbr$
 with real analytic inverse $s\mapsto - \mu_k(x-sa)$ where $\mu_k$ is a $b$-eigenvalue.
 }
\medskip

The proof of this result uses only the classical elementary fact that each point on an
algebraic curve admits a local uniformizing parameter.  \Ga's theory,
together with a recent convexity result of Bauschke, G\"uler, Lewis and Sendov [BGLS],
is then presented as a series of elementary consequences of this theorem.  This is done 
in  Sections 2 and 3.  Another recent result, due to Gurvits [Gur$_{1,2}$], is presented in Appendix B as
an improvement of a basic inequality of \Ga. The parallels in the proofs are emphasized.

\medskip
\noindent
{\bf Historical Note:}  In 1958 P. Lax  [L] conjectured that for homogeneous polynomials
$p(x,y,z)$   in three variables which are hyperbolic with respect to a direction $a$, say 
$a=(0,0,1)$,  there exist symmetric matrices $A, B$ such that 
$
p(x,y,z) = \det(xA+yB+zI)
$.
 This conjecture was essentially 
established by  Vinnikov [V] in 1993 (see [LPR]). Armed with this highly non-trivial result, one can quickly establish the full \Ga\ theory, since all arguments only involve three variables at  a time.
However, the presentation here emphasizes quite elementary  proofs.
\medskip

The second objective of this paper is to investigate the relationship of G\aa rding theory to 
the theory of  fully nonlinear second-order (degenerate  elliptic) equations, first discussed by 
Cafferelli, Nirenberg and Spruck  [CNS] for O$_n$-invariant equations.
This   discussion can be easily  ignored by readers interested only the first  part of the paper
on \Ga\ theory.
Here the basic vector space in question is the space $V=\Symn$ of real symmetric $n\times n$ matrices.
This is because a second-order equation can be viewed as a subset $F\ss \Symn$.
The {\sl solutions} or {\sl $F$-harmonic functions} are functions $u$ satisfying (at least in the $C^2$-case)
$$
(D^2 u)(x) \ \in\ \partial F.
$$
In [HL$_1$] the Dirichlet problem is studied and solved via the Perron method utilizing
the  notion of an $F$-subharmonic function. 
The only condition the closed subset $F$ must satisfy in this theory is the positivity
condition
$$
F+ P \ \ss\ F \qquad{\rm for\ each\ \ } P>0 \ \ {\rm (positive\ definite).}
\eqno{(1.4)}
$$
In this case the closed set $F\ss\Symn$ is referred to as a {\sl subequation}.
One of the main points of this paper is to show how \Ga's theory determines a
vast array of interesting and natural subequations.  His monotonicity result
(1.3) is exactly what is needed to establish the crucial positivity condition (1.4).

The first and most basic example is $F=\overline\G$, the \Ga \ cone itself.  Here the positivity
condition (1.4) must be assumed, not deduced. However, requiring (1.4) is equivalent
to requiring that the polynomial $p$ be hyperbolic in each positive definite direction $P>0$
in $\Symn$. These polynomials will be called {\sl \DG polynomials} for the purposes of this paper.

The main new construction of subequations in this paper can be described as follows. 
Suppose that $p$ is a \DGP\  of degree $m$, and choose a hyperbolic direction
$a\in\G$, for example any $P>0$ or, more specifically, the identity $I\in\Symn$.
Then the eigenvalue map $\l:\Symn\to \bbr^m$ is defined, as above, modulo coordinate permutations
on $\bbr^m$.  Consequently each closed symmetric subset $\E\ss \rn$ determines a second-order
partial differential equation (at least for $C^2$-functions $u$) by requiring
$$
\l\left( (D^2 u)(x)\right) \ \in\ \partial \E.
\eqno{(1.5)}
$$

One main result (Theorem 5.19) says that if a closed symmetric subset 
$E\ss\bbr^m$ satisfies the {\sl universal positivity condition}:
$$
\E + \l\ \ss\ \E \fa \l \in\bbr^m \ {\rm with}\ \l >0 \ \ {\rm (i.e.,\ each\ } \ \l_k>0),
\eqno{(1.6)}
$$
then the subset $F_\E$ of $\Symn$ defined by 
$$
F_\E \ \equiv\ \{ A\ \in\ \Symn  :  \l(A)\ \in\ \E\}
\eqno{(1.7)}
$$
is a subequation, that is, $F_\E$ satisfies the positivity condition (1.4).
The closed subsets $\E$ satisfying (1.6) can be regarded as {\sl universal subequations
on $\bbr^m$} since they induce a subequation on $\rn$ for each \DGP\ of degree $m$ on $\Symn$
(and each hyperbolic direction). For each such subequation on $\rn$, Theorem 5.20
concludes that the Dirichlet problem has a unique continuous solution for all continuous
boundary functions on any smooth domain $\O$ in $\rn$ whose boundary is suitably
``pseudo''-convex.

 This applies in particular to each of the branches  of the equation
 $$
 p((D^2 u) \ =\ 0,
 $$
where the  {\bf  $k$th branch} $F_k$ 
is defined by requiring that 
$$
\lu_k( D^2 u) \ \geq \ 0
$$
where $\lu_k$ is the $k$th ordered eigenvalue function of $p$.
These branches $F_k$ are independent of the choice of hyperbolic direction
(Theorem 2.12).   The universal $k$th branch $\E_k\ss \bbr^m$ is defined by
$$
\E_k \ \equiv\ \{\l\in\bbr^m : {\rm at\ least\ } m-k+1 \ {\rm  of\  the \ } \l_j\ {\rm are\ } \geq0 \} \ =\ \{\l\in\bbr^m: \lu_k\geq0\}  
$$
  In Appendix A we give an algebraic description of each of these  branches in terms of polynomial
inequalities.

We shall return to a discussion of examples of universal eigenvalue subequations $\E\ss\bbr^m$, but first we 
present examples of \DGP s.
Three basic irreducible classical examples are the three determinants:
$\det_\bbr$ on $\Symn$, $\det_\bbc$ on $\Sym(\bbr^{2n})$ and 
$\det_\bbh$ on $\Sym(\bbr^{4n})$.
They are discussed in Section 5.
A new   irreducible example of degree $2^n$ on  $\Sym(\bbr^{2n})$,  which yields a notion
of Lagrangian harmonicity, is also discussed there. Our theory  enables us to treat all branches of this new equation
even though an explicit rendering of them is quite complicated.

Several methods of constructing new \DGP s from a given \DGP\  are outlined in Section 5 as well.
One  method is to take a directional derivative in one of the directions in the \Ga\ cone.
Equivalently, these polynomials are the elementary symmetric functions in the eigenvalue functions.
A  second method constructs a \DGP\  whose eigenvalue functions are the ${m\choose p}$ possible $p$-fold sums of the given eigenvalue functions.  A third method constructs (for each $\d>0$) a new \DGP\ 
(of the same degree) whose open \Ga\ cone is uniformly elliptic, i.e., contains the closed set of all
$A\geq0$ in $\Symn$.
Since implementing these methods is not a commutative process, taken together they yield a huge collection
of \DGP s.

The second and third methods are special cases of a more general method.  If $Q$
is a symmetric polynomial on $\bbr^m$ of degree $N$ which is hyperbolic in all directions $\l>0$,
and $p$ is a \DGP\  on $\Symn$, then $q(A) = Q(\l(A))$ is a \DGP\  of degree $N$ on $\Symn$.

Now we return to the other ingredient, the universal subequations $\E$.
A classification or ``structure'' theorem is given in Proposition 5.25 in terms of Lipschitz graphs.
However,  there are examples which are obscured by this structure theorem, such as the 
universal Special Lagrangian subequation $E_c\ss\bbr^m$ defined by requiring
that 
$$
\sum_{k=1}^m \arctan \l_k \ \geq \ c.
$$
The structure theorem says that every  universal subequation $E$ has  boundary which is 
a graph  
$$
\partial E \ =\  \{ f(\l) e + \l : \ \l\cdot e=0\}
$$
over the hyperplane $\{e\cdot \l=0\}$ where $e=(1,...,1)$. $E$ is then the region above the 
graph in this picture.  The graphing function $f$ must be Lipschitz with Lipschitz constant
one using a natural norm on $\{e\cdot \l=0\}$ .

The Dirichlet problem is solved in Theorem 5.20 without assuming that the subequation
$F_E$ defined by (1.7) is either convex or uniformly elliptic, but our construction includes many
such examples because, if $E$ is convex, then $F_E$ is convex, and if $E$
is uniformly ellpitic, then so is $F_E$ (see Theorem 5.19 (a) and (b)).

The general G\aa rding theory has an important application to the analogous Dirichlet problem
on riemannian manifolds.  In a recent paper [HL$_2$] the authors have carried over
the basic results of  [HL$_1$] from constant coefficient purely second-order subequations
 to equations defined in a quite general setting on
riemannian manifolds, and manifolds with reduced structure group --
such as almost complex hermitian manifolds.  A fundamental new feature of this
theory is the introduction of a {\sl monotonicity cone}  $M_F$ for a given subequation
$F$.  Existence and uniqueness results are established under the (essentially necessary)
hypothesis that there exist some strictly $M_F$-subharmonic function defined on a 
neighborhood of the domain in question. The \Ga \ monotonicity result (1.3) exactly provides such
monotonicity cones.  For a given $I$-hyperbolic polynomial $p$ on $\Symn$ the 
\Ga \ cone $\G$ is a monotonicity cone for each branch of the equation $\{p(D^2 u) =0\}$.
More generally, the subequations $F_E$ on $\bbr^n$, determined by $p$ of degree $m$
and   a universal subequation $E$ on $\bbr^m$, are not confined to domains in $\rn$, 
but extend to all parallelizable riemannian $n$-manifolds, and when $p$ has a non-trivial
invariance group $G$, they extend to all riemannian $n$-manifolds with  structure group reduced (topologically) to $G$.
Details are discussed in Chapter 6.

\vfill\eject


\centerline{\headfont \ 2.\  G\aa rding's Theory of Hyperbolic Polynomials}
\medskip

This section gives an essentially self-contained account of the G\aa rding theory of hyperbolic
polynomials.  Some aspects of the presentation are new -- for example Theorem 2.7. Its
 proof relies only on   the classical elementary 
 fact that each point on an algebraic curve
admits  local uniformizing parameters (cf. (2.12)). This Theorem 2.7  together with
Remark 2.3 imply all the basic results in [G$_2$].
 
 \bigskip
\centerline{\headfont   Preliminaries}
\medskip

Suppose that $p$ is a homogeneous polynomial of degree $m$ on a complex vector space $V_\bbc$. 
Given points
$a,x\in V_\bbc$ with $p(a)\neq0$, the one-variable polynomial $s\mapsto p(sa+x)$ factors as 
$$
p(sa+x)\ =\ p(a)\prod_{k=1}^m \left( s + \l_k(x) \right) \qquad {\rm and}\ \ p(x)\ = \ p(a)\prod_{k=1}^m \l_k(x),
\eqno{(2.1)}
$$
where, by definition, $\l_1(x),...,\l_m(x)$ are the {\bf $a$-eigenvalues of $x$} and $r_k(x)\equiv -\l_k(x)$ are 
the {\bf  $a$-roots of $x$}.
Of course $\l(x) \equiv (\l_1(x),...,\l_m(x))$  is defined only  modulo the 
permutation group $\pi_m$.  When necessary, the dependence of the eigenvalues
of $x$ on the choice of $a$ will be denoted by $\l^k_a(x) \equiv  \l_k(x)$
and $\l_a(x) = (\l^1_a(x),...,\l_a^m(x))$, and similarly for the $a$-roots of $x$.

Everything proven in this section uses one or more of the elementary properties discussed next -- 
and very  little  else.

\medskip
\noindent
{\bf Elementary Properties:}  Suppose $p(a)\neq 0$.  It follows easily from (2.1) that for all 
$t\in \bbc$ and $x\in V_\bbc$:
\medskip

\qquad (1) \ \ $\l^k_a(t x)\ =\  t  \l^k_a(x)  \  {\rm mod  } \ \pi_m$,

\medskip

\qquad (2)  \ \ $\l^k_a(t a+x)\ =\  t+\l^k_a(x)  \  {\rm mod  } \ \pi_m$, \ and so
\  (2a)  \ $\l^k_a(a)\ =\  1$ \ for  all $k$.
\medskip

Sometimes it is convenient to set $e=(1,...,1)$ and restate (2) in the form

\medskip

\qquad (2)$'$  \ \ $\l_a(t a+x)\ =\  te+\l_a(x)  \  {\rm mod  } \ \pi_m$.
\medskip

 Using (2.1) with $x=tb$, $t\in \bbc$, $b\in V_\bbc$ and   applying (1) yields
 $$
 p(sa+tb)\ =\ p(a)\prod_{k=1}^m \left( s + t \l_k(b) \right)  
\eqno{(2.1)'}
$$
 
Note that (2.1)$'$ gives an explicit factorization of $p$ restricted to 
 $\span\{a,b\}$. This will be useful in establishing certain important  properties.
 \medskip
 
 We shall use the notation
 $$
 p_x'(y)\ \equiv\ {d\over dt}p(y+tx)\bigr|_{t=0}
 $$
 for the directional derivative of $p$ at $y$ in the direction $x$.  Setting
 $s=1$ in (2.1)$'$ and then taking the logarithmic derivative at $t=0$ proves that
 $$
  (3) \ \ \ p(a+tx) \ =\ p(a) \prod_{k=1}^r \left( 1+t \l^k_a(x)  \right)
 \and
 { p_x'(a) \over p(a) }  \ =\ \sum_{k=1}^r \l_a^k(x). \ \ 
$$
Here $r=r(x)$, called the {\bf rank of $x$}, is the number of non-zero $a$-eigenvalues, 
 which are here listed first.  The {\bf nullity of $x$} is the number $n(x) = m-r(x)$ of 
 $a$-eigenvalues which are zero. The elementary formula (3) is particularly  useful for studying degenerate
 cases where the rank of $x$ is less than $m$, i.e., where $p(x)=0$.

  As a consequence of the second part of Property (3) we have
  \medskip

\qquad (4) \ \  
 The $a$-trace of $x$, trace$_a(x) \ =\ \sum_{k=1}^r \l_a^k(x)$ is a linear functional in $x$.
\medskip

A  property, not used here and left as an exercise,  is that for $b\in V_\bbc$ with $p(b)\neq 0$,
$$
 \l^k_b(a) \ =\  {1\over \l^k_a(b)} \  {\rm mod  } \ \pi_m.
$$

The standard algebraic fact that the roots of a polynomial depend continuously on the coefficients
(if the degree does not drop) will be used in several guises.  For instance,

\medskip

\qquad (5) \ \ The $\l^k_a(x) \  {\rm mod  } \ \pi_m$ are continuous in $x\in V_\bbc$.
\medskip


 \bigskip
\centerline{\headfont   Hyperbolic Polynomials}
\medskip

Now suppose that  $p$ is a homogeneous real polynomial of degree $m$ on  a real vector space $V$.

\Def{2.1}  Given $a\in V$, the polynomial $p$ is {\bf  $a$-hyperbolic} if $p(a)>0$ and the 
eigenvalues $\l_k(x)$, $k=1,...,m$ are real for all $x\in V$.  That is, the 
 polynomial $t\mapsto p(ta +x)$ has exactly $m$ real roots for each $x\in  V$.
 The {\bf G\aa rding cone} $\G\ss V$ (also denoted by $\G_a$) is defined to be
  the set of $x$ such that $\l_k(x)>0$ for $k=1,...,m$.   

 \medskip
 
 Of course, the elementary properties (1) through (5) hold for $p$ on $V$.

  
\vskip .3in 
\centerline{\headfont A Basic  Inequality}

\medskip

Here we prove an elementary version of \Ga's inequality.

\Lemma{2.2}  {\sl
Suppose $p$ is $a$-hyperbolic of degree $m$.  If $x\in\G_a$, then}
$$
\left({p(x)\over p(a)}  \right)^{1\over m} \ \leq \ {1\over m} {p_x'(a) \over p(a)},
\qquad {\rm with\ equality\  } \iff\ \l_1(x)=\cdots \l_m(x)
\eqno{(2.2)}
$$
\pf
Each of the $a$-eigenvalues $\l_k$ of $x$ is $>0$. Using the second part of (2.1) and of property (3),
one sees that (2.2) is the classical inequality between the geometric and the arithmetic mean:
$$
\left(   \prod_{k=1}^m\l_k \right)^{1\over m} \ \leq \ {1\over m}  \sum_{k=1}^m\l_k,
\qquad {\rm with\ equality\  } \iff\ \l_1(x)=\cdots \l_m(x)  \qquad \mathqed
\eqno{(2.2)'}
$$

\Remark{2.3}  The inequality (2.2) has a geometric formulation. 
 The condition that the graph of a function $f(x)$ lies below its
 tangent plane at a point $(a,f(a))$ can be written as:
$$
f(x) \ \leq \ f(a) + f_{x-a}'(a) \ =\ f(a) -f_a'(a) + f_x'(a) \qquad {\rm for\ } x \ {\rm near\ } a.
\eqno{(2.3)}
$$
If $f$ is homogeneous of degree one, recall that $f_a'(a) = f(a)$, so that condition (2.3) simplifies to
$$
f(x) \ \leq \  f_x'(a) \qquad {\rm for\ } x \ {\rm near\ } a.
\eqno{(2.3)'}
$$
Take $f(x) = p(x)^{1\over m}$ with $x\in\G_a$.  Then $f_x'(a) = {1\over m} p(a)^{{1\over m}-1} p_x'(a)$ and we
see that (2.3)$'$ and (2.2) are the same inequality.  This proves that
$$
{\rm If\ } p \ {\rm is\ } b \ {\rm hyperbolic\ for\ all\ } b\in\G_a,\ {\rm then\ }
p(x)^{1\over m} \ {\rm is \ a\ concave\ function\  on\ }\G_a.
\eqno{(2.4)}
$$
In particular, in this case $\G_a$ is convex, since $-p(x)^{1\over m}$ is convex and negative on $\G_a$ 
with boundary values zero.

  
\vskip .3in 
\centerline{\headfont Ordered Eigenvalues}

\medskip

Suppose $p$ is an $a$-hyperbolic polynomial. 
Given a vector $\l\in\bbr^m$, let $\lu\in \bbr^m$ denote the non-decreasing reordering of $\l$, i.e., 
$
\lu_1 \ \leq\ \lu_2 \ \leq\ \cdots \ \leq\ \lu_m.
$
 The {\bf ordered eigenvalue functions}
are then defined to be $\lu(x) = (\lu_1(x),...,\lu_m(x))$.
Note that  $\l_{\rm min}(x) = \lu_1(x)$ and 
$\l_{\rm max}(x) = \lu_m(x)$.
It is useful to rewrite some of the properties listed above
in terms of the ordered  $a$-eigenvalues.

\medskip\noindent
{\bf Elementary Properties:}   For each $k=1,...,m$, we have
\medskip

(1a)\ \ $\lu_k(tx)\ =\ t\lu_k(x)$,  for $t\geq 0$,  
\medskip

(1b)\ \  $\lu_k(-x) = -\lu_{m-k+1}(x)$, and in particular, $\l_{\rm max}(x) = - \l_{\rm min}(-x)$
\medskip

(2)\ \ $\lu_k(ta+x)\ =\ t+\lu_k(x)$ for all $t\in \bbr$,
\medskip


(5)\ \   $\lu_k(x)$ is  continuous.
\bigskip

\Def{2.4}  The sets $$F_k \equiv \{ \lu_k \geq 0\}$$ (rather than $\partial F_k$) will be referred to as the
{\bf branches} of $\{p=0\}$.  The {\bf principal} or {\bf smallest branch} 
$F_1 = \{\l_{\rm min}\geq 0\}$ is of particular interest.  
\medskip

Note that each branch $F_k$ is a cone with vertex at the origin by  (1a).

\Remark{2.5}
Given a continuous function $g:\bbr^{N-1}\to \bbr$, we define {\bf the associated
graphing function} $G:\bbr^N\to\bbr$  by $G(t,x) = t-g(x)$.  The closed set $F\equiv \{G\geq 0\}$
above the graph of $g$ has interior $\Int F = \{G > 0\}$ and boundary
$\partial F = \{G = 0\}$  equal to the graph of $g$.  Moreover, $\Int F$ is a connected open set 
with $F= \overline{\Int F}$.
Finally note:
$$
F \ \ {\rm  is \ convex\ } \iff \  g  \ \ {\rm  is \ convex\ } \iff \ G \ \ {\rm  is \ concave\ }
\eqno{(2.5)}
$$

These obvious facts can be applied as follows.  
Suppose that $p$ is an $a$-hyperbolic polynomial and consider
the ordered eigenvalues.  Take $\bbr^{N-1}$ to be a hyperplane $W\ss V$
transverse to $a$.
\medskip
\centerline
{
If  $g$  is  the restriction of $-\lu_k$ to $W$, \ then  $G=\lu_k$
}
\medskip
\noindent
 since for
 $x\in W$, we have
$
G(ta+x)\ =\ t-g(x) \ =\ t+\lu_k(x)\ =\ \lu_k(ta+x)
$
where the last equality follows from Property (2). Consequently, as a special case of 
Remark 2.5 we have:

\Prop{2.6}  {\sl  Suppose that $p$ is $a$-hyperbolic. 
Then $\partial F_k=\{\lu_k=0\}$, and $\Int F_k = \{\lu_k>0\}$.  Furthermore, $\Int F_k$  is a connected set satisfying
$F_k=\overline{\Int F_k}$.  
In particular
$$
\G\ =\ \Int F_1,  
\qquad 
F_1\ =\ \overline \G, \qquad {\rm and}
$$
$$
\G \ \ {\rm  is \ convex\ } \iff \  \l_{\rm min}  \ \ {\rm  is \ concave\ } \iff \  \l_{\rm max}  \ \ {\rm  is \ convex\ }
\eqno{(2.6)}
$$
}
\medskip

Thus the G\aa rding cone is the interior of the principal branch.  The last equivalence in (2.6)
follows since $\l_{\rm max} (x) = - \l_{\rm min} (-x)$ by   (1b).

\Cor{2.7} {\sl The G\aa rding cone $\G_a$ is the connected component of $\{p\neq 0\}$ which contains $a$.}

\pf  Let $\G$ denote the connected component of $\{p\neq 0\}$ containing $a$.
By Property (2a) we have $a\in\G$, and since $\G$ is open and connected, $\G_a\ss \G$.
Since the $a$-eigenvalues $\l_k(x)$ never vanish on $\G$ and equal 1 at $a\in \G$, they are strictly
positive on $\G$, i.e., $\G \ss\G_a$.\qed
 
\Remark{2.8}  Setting $k^* = m-k+1$, we have
\medskip
\centerline
{
$x \in \Int F_k \quad\iff\quad x$ has at least $k^*$ strictly positive eigenvalues.
}
\medskip

  
\vskip .3in 
\centerline{\headfont A Real Analytic Arrangement of the Restricted Eigenvalues}

\medskip

Throughout this subsection we assume that $p$ is an $a$-hyperbolic polynomial of degree $m$ on $V$.
The eigenvalues $\l_1(x),...,\l_m(x)$ of $x\in V$ are well defined mod $\pi_m$.   The continuity of 
$\l_k$ as a function of $x$ (Property (5)) cannot be improved to differentiability.  However, by restricting
the $\l_k$ to the line through  $x$ in a direction $b\in \G$, the one-variable functions 
$
 \l_1(x+tb), \ ...\  ,  \l_m(x+tb),
$
of $t$ defined mod $\pi_m$, can be arranged so that each function  is a
strictly increasing bi-real-analytic mapping
of $\bbr$ onto $\bbr$.  
The functions
$$
\l_1(x-tb),\ \  ... \ \ , \l_m(x-tb)
$$
of $t$  defined mod $\pi_m$ will be referred to as the
{\bf restricted eigenvalue functions}.

\Theorem{2.9}  {\sl Suppose that $p$ is an $a$-hyperbolic polynomial of degree $m$ on $V$
with G\aa rding cone $\G$.  Then $p$ is also $b$-hyperbolic for each  $b\in \G$, and $\G_b=\G_a$.  
Moreover, for each fixed $x\in V$ the restricted eigenvalue functions of $t$
$$
\l^1_a(x-tb),\ \  ... \ \ , \l^m_a(x-tb)
\eqno{(2.7)}
$$
can be arranged so that
each one  is a  real analytic, strictly decreasing  function 
from $\bbr$ onto $\bbr$, and inverses are given by   
$$
\l^1_b(x-sa),\ \  ... \ \ , \l^m_b(x-sa).
\eqno{(2.8)}
$$
Consequently,}
$$
{ d \over dt} \l^k_a(x+tb) \ >\ 0 \fa t\in \bbr \ \ {\rm and\ \ } k=1,...,m.
\eqno{(2.9)}
$$

\medskip
\noindent
{\bf Remark.}  Note that by combining Corollary 2.7 and Theorem 2.9 we have that
for each connected component $\G$ of the set $\{p\neq 0\}$, \medskip
\centerline
{\sl
$p$ is $a$-hyperbolic for some $a\in\G$ \quad$\iff$\quad $p$ is $a$-hyperbolic for all $a\in\G$.
}
\medskip\noindent
Consequently it is reasonable to use the phrase ``$p$ is $\G$-hyperbolic'' where 
$\G$ is a connected component of $\{p\neq0\}$.\medskip

 Theorem 2.9 can be deduced from a local result, which holds without assuming $b\in \G$.
 
 \Lemma{2.10} {\sl 
 Suppose $x, b \in V$. Then near any point $t_0\in \bbr$
 the restricted root functions
 $$
 r^1_a(tb+x),\ ... , r^m_a(tb+x)
 \eqno{(2.10)}
$$
 can be arranged to be real-analytic in $t$.
 Therefore, by uniqueness of real analytic continuation the local arrangements
 (2.10) extend to a global arrangement of the restricted roots, which we label as }
$$
s_k(t)  \ =\ r^k_a(tb+x) \qquad k=1,...,m.
\eqno{(2.11)}
$$
 \pf
 Let $s_1$ denote one of the root functions evaluated at $t_0$.
 That is, $(s_1, t_0)$ is a point on the complex algebraic curve
$$
C\ \equiv\ \{(s,t) : p(sa+tb+x)=0\}
$$
 in $\bbc^2$.  Choose a local irreducible component $C_1$ of $C$ at 
 $(s_1, t_0)$, along with a local uniformizing parameter $z$.
  This means that 
we have holomorphic functions 
$$
s(z) \ =\ s_1 +z^pg(z) \and t(z) \ =\ t_0 + z^qh(z),\quad {\rm for} \ |z|<\d,
\eqno{(2.12)}
$$
with  $g(0)\neq 0$, $h(0)\neq 0$, such that $z \mapsto (s(z), t(z))$
is a homeomorphism of a neighborhood of $0\in\bbc$ onto a neighborhood  of
$(s_1,t_0)$ in $C_1$.  By extracting a $q$th root of $h(z)$ we may assume that $h(z) \equiv 1$, i.e.,
$$
t(z) \ =\ t_0 + z^q.
$$

It will suffice to show that $q=1$, for then
 $s(z)$ is a root of $p(sa+(t_0+z)b+x)$, and, as such, must be real if $z$ is real.
Setting
$$
r^1_a(tb+x)\ \equiv\ s(t-t_0)
$$
we obtain a local real analytic restricted root function.  Repeating  this construction 
for each of the $m$ real points $(s_1,t_0), ..., (s_n,t_0)\in C$ (counting multiplicities)
will complete the proof.

We now show that $q=1$.
Since $s(z)$ is a root of $s\mapsto p(sa+t(z)b+x)$ and $p$ is $a$-hyperbolic,
$$
z^q\ \ {\rm real} \quad\Rightarrow\quad  s(z)\ \ {\rm real}.
\eqno{(2.13)}
$$
In particular, $z$ real $ \quad\Rightarrow\quad $ $s(z)$ real, so that the power series
$$
s(z) = \sum_{k=0}^\infty a_k z^k
$$
for $s(z)$ has real coefficients $a_k$.  Set $\o \equiv e^{\pi i\over q}$ and $z=r\o$ with
$r$ real.  Then $\o^k$ is not real unless $k$ is a multiple of $q$.  Now 
$z^q$ is real and hence $\sum_{k=0}^\infty a_k \o^k r^k$ is real by (2.13).
It follows that $a_k=0$ unless $k$ is a multiple of $q$.   
Thus $s(z) = s_1+ \sum_{k=1}^\infty b_k z^{kq} = s_1+f(z^q)$.
Since $z\mapsto (f(z^q), z^q)) = (s(z)-s_1, t(z)-t_0)$ is one-to-one, this proves that $q=1$.
\qed

\medskip
\noindent
{\bf Proof of Theorem 2.9.} 
Now we will use the fact that $b\in\G$.
By Property (1)
$$
r^k_a(tb+x) \ =\  t r^k_a  \left (b+{x\over t}\right)\quad {\rm mod}\ \pi_m\qquad {\rm for\ }t\neq 0.
\eqno{(2.14)}
$$
For $|t| >>0$ sufficiently large,
$$
r^k_a  \left (b+{x\over t}\right)\ <\ 0.
\eqno{(2.15)}
$$
This follows from the continuity of the eigenvalues (5) and the fact that $r_a^k(b) < 0$, i.e., $b\in\G$.
Consequently,
$$
\lim_{t\to \pm \infty}  r^k_a(tb+x) \ =\ \mp\infty.
\eqno{(2.16)}
$$
In particular, none of the functions $s_k(t)$ are constant, and each maps onto all of $\bbr$.

Suppose that $\vf(t)$ is one of the function $s_k(t)$ which is repeated $\ell$-times.  
If $s_0$ is not a critical value of $\vf$, then 
$$
{\rm each \ point\ } t_j\in \vf^{-1}(s_0)  \ {\rm 
is\  a \ root\  of  } \ t\mapsto p(s_0a +tb +x) {\ \rm  of\  multiplicity } \ \ell. 
\eqno{(2.17)}
$$
 This follows from the 
fact that the functions $s-\vf(t)$ and $t-\vf^{-1}(s)$ vanish to first order on the same real curve
in $\bbr^2$ near $(s_0,t_0)$ and  therefore
differ by a smooth (in fact, real analytic) factor $\a(s,t)$ near $(s_0,t_0)$ with $\a(s_0,t_0)\neq0$.

Let   $S_0$ denote the set of points which are not   critical values for any of the functions $s_k(t)$.  
If $s_0\in S_0$, then $s_k^{-1}(s_0)$
can have at most one point $t_k$.  Otherwise by (2.17)  the polynomial $t\mapsto p(s_0a+tb+x)$ would
have more than $m$ roots.  Hence each $s_k$ is one-to-one outside its critical points.
This is enough to conclude that each $s_k$ is a homeomorphism from $\bbr$ onto $\bbr$.
Therefore, each $s_k(t)$ defined by (2.11) is a strictly  decreasing homeomorphism from $\bbr$ onto $\bbr$.

Moreover,  for $s_0\in S_0$ the points $t_1=s_1^{-1}(s_0),..., t_m=s_m^{-1}(s_0)$,
listed to multiplicity, provide $m$ (real) roots of $t\mapsto p(s_0a+tb+x)$.
This proves that
$$
r_b^k(sa+x)\ =\ s_k^{-1}(s)
\eqno{(2.18)}
$$
 for $s\in S_0$.  By continuity of the roots (2.18)
 remains true for critical values $s\in \bbr$.

In particular, the roots
$$
r^1_b(x) \ =\ s_1^{-1}(0),\ ...\ , r^m_b(x) \ =\ s_m^{-1}(0)
$$
of $p(tb+x)$ are all real, proving that the polynomial $p$ is $b$-hyperbolic.
By (2.16) the monotone functions $s_k(t) = r_a^k(tb+x)$ must be non-increasing homeomorphisms
of $\bbr$ onto $\bbr$.  Since $p$ is $b$-hyperbolic, the results established for $r_a$ apply
to $r_b$.  In particular, by (2.18) each of the inverses $s_k^{-1}$ is real analytic.  This proves 
that for some arrangement of the restricted roots
$$
s \ =\ r_a^k(x+tb) 
\and
t \ =\ r_b^k(x+ta) \ \ \ k=1,...,m
\eqno{(2.19)}
$$
are inverses of each other.  Also,
\medskip
\centerline
{
Each $r_a^k(x+tb)$, as a function of $t$, is a } 
\centerline
{\qquad\qquad\qquad\quad
strictly decreasing bi-real-analytic map from 
$\bbr$ onto $\bbr$. 
\qquad\qquad
(2.20)
}
\medskip
Finally note that (2.19) and (2.20) remain true with $r_a^k(x+tb)$ replaced by
$\l_a^k(x-tb) = -r_a^k(x-tb), \ k=1,...,m$. \qed

\Ex{2.11. (Comparing Real Analytic and Ordered Eigenvalues)}
The polynomial $p(x) = x_1^2-x_2^2-x_3^2$ on $\bbr^3$ is the prototype
of all hyperbolic polynomials.  It is $(1,0,0)$-hyperbolic with G\"arding cone equal
to the light cone $\G =  \{p(x)>0\} \cap \{x_1>0\}$. The ordered eigenvalue functions are:
$$
\l_{\rm min}(x) \ =\ x_1-\sqrt{x_2^2+x_3^2} \  \ \ \  \leq\ \ \ \ 
x_1+\sqrt{x_2^2+x_3^2} \ =\ \l_{\rm max}(x).
$$
Choose $b\in \G$ of the form $b_1> b_2 >0 = b_3$.  The ordered arrangement of the restricted
eigenvalue functions is
$$
\eqalign
{
\l_{\rm min}(x+tb) \ &=\ x_1 +tb_1-\sqrt{(x_2+tb_2)^2+x_3^2}  \cr
\l_{\rm max}(x+tb) \ &=\ x_1 +tb_1+\sqrt{(x_2+tb_2)^2+x_3^2}  \cr
}
$$
If $x_3\neq0$, then both of these functions are real analytic in $t$, and hence the 
real analytic arrangement agrees with the ordered arrangement.  If $x_3=0$, the ordered arrangement is:
$$
\l_{\rm min}(x_tb) \ =\ x_1 +tb_1 - |x_2 +tb_2| \  \ \ \  \leq\ \ \ \ 
\ x_1 +tb_1  +   |x_2 +tb_2|   \ =\ \l_{\rm max}(x+tb).
$$
while the real analytic arrangement is 
$$
\eqalign
{
\l^1(x+tb) \ &=\ x_1 +tb_1- (x_2+tb_2)  \cr
\l^2(x+tb) \ &=\ x_1 +tb_1+ (x_2+tb_2).  
}
$$
It is interesting to note the following asymptotics:
$$
\ \ \ \ \ \  {\rm for}\ \ t>>0,\qquad
\l_{\rm min}(x+tb) \ \sim\  t\l_{\rm min}(b) \and \l_{\rm max}(x+tb) \ \sim\  t\l_{\rm max}(b) 
$$
$$
{\rm while\ for\ \ } t<<0, \qquad 
\l_{\rm min}(x+tb) \ \sim\  t\l_{\rm max}(b) \and \l_{\rm max}(x+tb) \ \sim\  t\l_{\rm min}(b).
$$
which shows that (2.14) only holds mod $\pi_m$.
 \medskip

Since the $a$-eigenvalues are all real, the rank of $x\in V$ can be refined into the sum of the 
{\bf plus rank}, denoted $r^+(x)$, and the {\bf minus rank}, denoted $r^-(x)$.

\Theorem{2.12} {\sl
The quantities $r^\pm(x)$ and $n(x)$ are all independent 
of the direction $a\in \G$ chosen to compute the 
eigenvalues of $x$.  In particular, 
all of the branches $F_k\equiv \{\lu_k\geq0\}$ of $\{p=0\}$ are independent 
of the point $a\in\G$ chosen to compute the eigenvalues.}

\pf 
Choose $a,b\in\G$ and assume $r^+_a(x) = k^*$ (cf. Remark 2.8).
Let $\g_1,...,\g_m$ denote the $m$ curves in the $(t,s)$-plane
 defined by (2.19).  Exactly  $k^*$ 
 of these curves cross the negative $s$-axis.  By (2.20) this implies that at exactly  $k^*$ 
 of these curves cross the negative $t$-axis, i.e., $r^+_b(x) = k^*$.
The proof is similar for $r^-(x)$ and $n(x)$.
\qed

  
\vskip .3in 
\centerline{\headfont Convexity}

\medskip

\Cor{2.13} {\sl 
\medskip
{\rm (a)} \ \  The G\aa rding cone $\G$ is convex.  
\smallskip
{\rm (b)} \ \  $\l_{\rm min}(x)$ is  concave.

\smallskip
{\rm (c)} \ \  $\l_{\rm max}(x)$ is  convex.

\smallskip
{\rm (d)} \ \ $p(x)^{1\over m}$ is concave on $\G$.

\smallskip
{\rm (e)} \ \ $-\log \,p(x)$ is convex on $\G$.
}

\pf
The conditions (a), (b) and (c) are equivalent by Proposition 2.6.  
Condition (d) follows from (2.4) in Remark 2.3.
Since $-p(x)^{1\over m}$ is convex on $\G$, negative on $\G$, and zero on $\partial \G$,
condition (a) follows. Finally, recall that if $\vf$ is increasing and convex, then $\vf(f(x))$
is convex if $f$ is convex.  To prove (e) take $\vf(t) = -\log(-t)$ for $t<0$ and note that 
$\vf(-p(x)^{1\over m}) = -{1\over m} \log\, p(x)$. \qed

\medskip\noindent
{\bf Second Proof of (a).}
Suppose $x,y\in \G$, i.e., $\l_{\rm min}(x)>0$ and 
$\l_{\rm min}(y)>0$.  We must show that $\l_{\rm min}(sx+ty)>0$ for all $s,t>0$.
By Theorem 2.9,  $\G_x=\G_a$ so that  we can assume the eigenvalues are computed using
$a=x$.  Then by Properties (1a) and (2)
$$
\qquad\qquad\qquad\qquad\qquad\qquad
\l_{\rm min}(sx+ty)\ =\ s+t\l_{\rm min}(y)\ >\ 0. 
\qquad\qquad\qquad\qquad\qquad\qquad \vrule width5pt height5pt depth0pt
$$

\Remark{2.14}  For $x\in W$, a hyperplane transverse to $a$,  we may apply Remark 2.5 to 
$g(x) \equiv  -\l_{\rm min}(x)$, $x\in W$.  Since $g$ is convex, we prefer the following  norm notation:
$$
\quad
\|x\|^+ =\ -\l_{\rm min}(x) \and
\|x\|^- =\  \|-x\|^+    \ =\ \l_{\rm max}(x) \qquad {\rm for\ } x\in W
$$
where the last equality follows from Property (1b).
By Remark 2.5  we have  
$$
\eqalign
{
\partial \G \ &=\ \{\|x\|^+a+x : x\in W\} \ \ {\rm is\ the\  graph\ of\ }\ \|\bullet\|^+, \ \ {\rm and} \cr
-\partial \G \ &=\ \{- \|x\|^+a+x : x\in W\} \ \ {\rm is\ the\  graph\ of\ }\  - \|\bullet\|^- 
}
$$
over the hyperplane $W$.  The convexity of $\G$ and 1(a) imply that 
$$
\|x+y\|^\pm \ \leq\ \|x\|^\pm+\|y\|^\pm
\and
\|tx\|^\pm \ =\ t\|x\|^\pm \fa t>0.
$$

  
\vskip .3in 
\centerline{\headfont Monotonicity}

\medskip

\Theorem{2.15}  {\sl The ordered eigenvalues are strictly $\G$-monotone, that is}
$$
\lu_k(x+b) \ > \ \lu_k(x) \fa x\in  V, b\in \G.
$$

\pf
Under the real analytic arrangement in Theorem 2.9   each of the restricted  eigenvalue functions 
$\l_k(tb+x)$ is strictly increasing in $t$.  Hence the ordered eigenfunctions

$$
\lu_1(tb+x)\ \leq\ \cdots\ \leq\ \lu_m(tb+x)
$$
are also strictly increasing.  Now compare $t=0$ with $t=1$.\qed

\Def{2.16}  We say that a subset $F\ss V$ is {\bf $\G$-monotone} if
$$
F+\G\ \ss\ F \quad {\rm or\ equivalently\ \ \ } F+\G\ \ss\ \Int F.
\eqno{(2.21)}
$$
The equivalence follows since $F+\G$ is open.  

\Cor{2.17}  {\sl
Each branch $F_k \equiv \{\lu_k \geq 0\}$ of $\{p=0\}$ is $\G$-monotone.
}

  
\vskip .3in 
\centerline{\headfont Constructing Hyperbolic Polynomials}

\medskip

In this subsection we describe some of the  important ways of constructing new hyperbolic polynomials
from a given hyperbolic polynomial.

\medskip
\noindent
{\bf I. Factors and Products.}  Suppose
 $p(x)$ and $q(x)$ be real homogeneous polynomials  with
 $p(a), q(a)>0$ for fixed $a\in V$.  Then the following is obvious.
 
\Prop{2.18} \medskip
\centerline{\sl $p(x)q(x)$ is $a$-hyperbolic \ \ $\iff$\ \ $p(x)$ and $q(x)$ are each $a$-hyperbolic.}
\medskip
\centerline{\sl     in which case \ \  $\G(pq) = \G(p)\cap \G(q)$.}
\medskip

\noindent
{\bf II. Restriction.} Suppose that $W$ is a vector subspace of $V$ and that $a\in W$.
Then the following is also obvious.

\Prop{2.19} \medskip
\centerline{\sl  If $p(x)$ is $a$-hyperbolic, then
$p\bigr|_W$ is $b$-hyperbolic for each $b\in \G(p)\cap W$
}
\medskip
\centerline{\sl     in which case \ \  $\G\left(p\bigr|_W\right) = \G(p)\cap W$.}
\medskip

\noindent
{\bf III. Closure.} 
Suppose $\{p_j\}$ is a sequence of $a$-hyperbolic polynomials of degree $m$
converging to a polynomial $p$.

\Prop{2.20} {\sl
If $p(a)\neq 0$, then $p$ is $a$-hyperbolic.
}
\pf
Since $p(a)\neq0$, the $a$-eigenvalues of $x\in V$ are defined as complex 
numbers by (2.1).  By continuity of the eigenvalues mod $\pi_m$,
the limiting $a$-eigenvalues of $p$ must be real.  Obviously $p(a)>0$ and degree of $p$
equals $m$.\qed

\medskip

\noindent
{\bf IV. Derivatives.} 
Let $p'$ or $p_b'$ denote the  directional  derivative of  $p$  in the direction $b$, i.e.,
$$
p'(x) \ =\ {d\over ds} p(x+sb)\bigr|_{s=0}
\eqno{(2.22)}
$$
   
   \Prop{2.21} {\sl
   Suppose $p$ is $a$-hyperbolic of degree $>1$.
   If $b\in\G(p)$, then $p_b'$ is $a$-hyperbolic.
   That is, $\G(p)\ss \G(p_b')$.
   Moreover, for any $x$ the ordered  $b$-eigenvalues for $p_b'$ 
   are interspersed between the ordered  $b$-eigenvalues for $p$.
   }
    \pf
    Note that by homogeneity $p_b'(b) = mp(b) >0$.
   Also $p_b'(tb+x) = {d\over dt} p(tb+x)$, and therefore by Rolle's Theorem
    the $b$-eigenvalues of $x$ for $p'$ are interspersed between the    
   $b$-eigenvalues of $x$ for $p$.  Namely, we have
   $$
\l_1(x) \leq \l_1'(x) \leq \l_2(x) \leq \cdots  \leq \l_{m-1}(x) \leq \l_{m-1}'(x) \leq \l_{m}(x).
\eqno{(2.23)}
$$
More sharply stated, 
\medskip
\centerline
{\qquad\qquad (a)\qquad
 $\l_k(x) < \l_{k+1}(x)\quad \Rightarrow\quad  \l_k(x) < \l_k'(x) < \l_{k+1}(x)$, \quad  and
\qquad\qquad\qquad (2.23)$'$
 }
\centerline
{ (b)\qquad  $\l_k(x)$ has multiplicity $\ell \quad\Rightarrow\quad
   \l_k'(x)$  has  multiplicity $\ell-1$. \qquad\qquad
   }
   \medskip
   \noindent
     Thus the $b$-eigenvalues of $x$ for $p_b'$
   are all real, proving that $p_b'$ is $b$-hyperbolic.
   In particular, by (2.23), if $x\in \G(p)$ (i.e., $\l_k(x)>0$, $1\leq k\leq m$), 
   then $\l_k'(x) >0$, for $1\leq k\leq m-1$ (i.e., $x\in \G(p')$).
   Thus $\G(p)\ss\G(p_b')$ proving that $p_b'$ is $a$-hyperbolic.
   \qed
   \medskip
   
   By induction we have
   \Cor{2.22} {\sl
   Suppose that $p$ is $\G$-hyperbolic and $b_1,...,b_k \in \G$.
   Then the $k$-fold directional derivative $p_{b_1,...,b_k}^{(k)}$ in the directions 
   $b_1,...,b_k$ is $\G$-hyperbolic, i.e.,  $\G(p)\ss \G(p^{(k)})$.
   }
   \medskip

\noindent
{\bf V. Elementary Symmetric Functions.} 

Let $\s_k$ denote the $k$th elementary symmetric function of $m$ variables.
If $p$ is a homogeneous polynomial of degree $m$ with $p(a)\neq0$, 
then   by (2.1) 
$$
p(sa+x) \ =\ p(a)\left[ s^m + \s_1 \left( \l(x)\right) s^{m-1} +\cdots+
\s_m \left(\l(x)\right)    \right].
\eqno{(2.24)}
$$
Taking directional derivatives in the direction $a$ we have
$$
\eqalign
{
p^{(k)}(sa+x) \ 
&=\ p(a)\left[ \a_0 s^{m-k} + \a_1 \s_1(\l(x))s^{m-k-1} +\cdots + \a_{m-k}\s_{m-k}(\l(x))      \right]
}
\eqno{(2.25)}
$$
with $\a_j = k!{m-j \choose k}$.
Setting $s=0$
 proves that the $k$th directional derivative of $p$ in the direction $a$
is (up to a constant) the same as the $(m-k)$th elementary symmetric function of the 
$a$-eigenvalue functions, i.e., 
$$
p^{(k)}(x)  \ =\  k! \, p(a) \s_{m-k}(\l(x)).
\eqno{(2.26)}
$$
In particular, this proves that
$$
\s_j(x) \ \equiv \ \s_j \left(\l(x)\right)
\ \ {\rm defines \ a\  homogeneous\  polynomial \ of\  degree\ \ }j
\eqno{(2.27)}
$$
and that its derivative in the direction $a$ is $(m-j+1)$ times $ \s_{j-1}(\l(x))$, i.e., 
$$
\s_j'(x) \ =\ (m-j+1) \s_{j-1}(x)
\eqno{(2.28)}
$$
By Corollary 2.22 this proves 
\Cor{2.23} {\sl
Suppose  $p$ is $a$-hyperbolic. Then $\s_k(x) = \s_k(\l(x))$ defines an $a$-hyperbolic polynomial
of degree $k$, and $\G(\s_{k+1})\ss \G(\s_k)$ for $k=1,...,m$.
}
\medskip

Let $\l^{(k)}_1(x), ..., \l^{(k)}_{m-k}(x)$ mod $\pi_{m-k}$ denote the $a$-eigenvalue functions for
the $a$-hyperbolic polynomial $p^{(k)}(x)$.  The equation (2.25) for  $p^{(k)}(sa+x)$ 
proves that the $j$th elementary symmetric functions
$$
\s_j(\l^{(k)}(x)) \ \ {\rm and\ \ } \s_j(\l(x)) \ \ {\rm are \  equal\ modulo\  a\  positive\  scale}.
\eqno{(2.29)}
$$


\medskip
\noindent
{\bf VI. Hyperbolic Polynomials Defined Universally.}  
Set $e=(1,...,1) \in \bbr^m$.
Suppose that $Q$ is a symmetric homogeneous polynomial of degree $N$ on $\bbr^m$.
If $p$ is a homogeneous polynomial of degree $m$ on $V$ with $a$-eigenvalue
functions $\l(x) = (\l_1(x),...,\l_m(x))$, then 
$$
q(x) \ =\ Q(\l(x))
\eqno{(2.30)}
$$
defines a homogeneous polynomial of degree $N$ on $V$
(since $Q$ is a polynomial in the $\s_\ell$'s and each $\s_\ell(\l(x))$ is  a polynomial).
If $Q$ is $e$-hyperbolic on $\bbr^m$, let $\L_1(\l),...,\L_N(\l)$ denote the $e$-eigenvalue functions for $Q$.

\Prop{2.24}  {\sl  Suppose $Q$ is a symmetric $e$-hyperbolic polynomial of degree $N$ on $\bbr^m$.  
Then $Q$ universally determines an $a$-hyperbolic  polynomial $q(x)$ defined by 
(2.30) on any vector space $V$ equipped with an $a$-hyperbolic polynomial $p(x)$ of degree $m$.
Moreover, the $a$-eigenvalue functions for $q$ are $\L_1(\l(x)),...,\L_N(\l(x))$.
}

\pf First note that $q(a) = Q(\l(a))=Q(e) > 0$.
By (2.1) 
$Q(te+\l) = Q(e)\prod_{j=1}^N(t+\L_j(\l))$.   Property (2)$'$:
$$
\l(ta+x) =\ te+ \l(x)
$$
implies that $q(ta+x) = Q(\l(ta+x))=Q(te+\l(x)) = Q(e)\prod_{j=1}^N(t+\L_j(\l(x)))$.
This proves that $\L_1(\l(x)),..., \L_N(\l(x))$, which are real, are the $a$-eigenvalues of the 
polynomial $q(x)$.\qed

\Remark{2.25}  Examples of $e$-hyperbolic polynomials $Q$ of degree $N$ on $\bbr^m$ abound.
Take any $a$-hyperbolic polynomial $R$ of degree $N$ on a vector space $W$.  Choose vectors
$b_1,...,b_m \in W$ with $b\equiv b_1+\cdots +b_m\in \G(p)$.
Then $Q(\l) = R(\l_1b_1+\cdots+\l_m b_m)$ is $e$-hyperbolic of degree $N$ on $\bbr^m$.
To use $Q$ in Proposition 2.24 it should first be symmetrized (thereby increasing its degree).

\medskip
Special cases of this construction V are of particular importance.
(The trivial case $Q(\l)=\l_1\cdots\l_m$ induces the polynomial $q=p$.)
 To begin exploring other cases note the following obvious

\medskip
\noindent
{\bf Fact 2.26.} {\sl  Given $w\in\bbr^m$ the linear polynomial $\ell(\l) \equiv w\cdot \l$ is
$e$-\hy if and only if $\ell(e) = w\cdot e = w_1+\cdots +w_m >0$, in which case the single 
$e$-eigenvalue function for $\ell$ is computed to be}
$$
\L(\l)\ =\ {w\cdot \l  \over w\cdot e}
\eqno{(2.31)}
$$

This polynomial is not symmetric.  Symmetrizing  $\ell(\l)$ yields 
$Q(\l) = \prod_{\s\in \pi_m} (\s w)\cdot \l$ and we obtain the following important
special case of Proposition 2.24.

\Prop{2.27}  {\sl
Suppose $p$ is $a$-hyperbolic of degree $m$ on $V$ with $a$-eigenvalue
functions $\l(x) = (\l_1(x),...,\l_m(x)) \ {\rm mod} \pi_m$.
Fix $w\in \bbr^m$ with $w\cdot e>0$. Then
$$
q(x) \ =\ \prod_{\s \in \pi_m} (\s w) \cdot \l(x)
$$
defines an $a$-hyperbolic polynomial on $V$ of degree $m!$ with $a$-eigenvalue functions}
$$
{(\s w)\cdot \l(x)  \over  w\cdot e}, \qquad \s\in\pi_m.
$$

The convexity of the largest eigenvalue function generalizes as follows.

\Cor {2.28. ([BGLS])} {\sl
Each non-decreaasing linear combination of the ordered eigenvalue functions is convex:
\medskip
\centerline{ $w\in\bbr^m_{\uparrow}$ \qquad $\Rightarrow$\qquad $w\cdot \lu(x)$ is convex on $V$.}
\medskip
\noindent
That is, for $x,y\in V$ and $0\leq s\leq 1$ one has }
$$
w\cdot \lu(sx+(1-s)y) \ \leq\ sw\cdot \lu(x) + (1-s)w\cdot \lu(y) \ =\ w\cdot\left(s\lu(x) + (1-s)\lu(y)   \right).
$$

\pf
Recall Property (4) that $e\cdot \lu(x) = \l_1(x)+\cdots+\l_m(x)$ is a linear function of $x$.
Therefore,
\smallskip
\centerline{
$w\cdot  \lu(x) \ \ {\rm is \ convex \ } \quad \iff \quad  (w+te)\cdot \lu(x)\ \ {\rm is \ convex, \ }$}
\smallskip
\noindent
and hence we may assume that $w\cdot e>0$. Then it is easy to see that 
 $w\cdot \lu(x)/ w\cdot e$ is the largest
eigenvalue function for $q(x)$ defined  in Proposition 2.27.
Therefore, $w\cdot \lu(x)$ is convex by  Corollary 2.13c.
 \qed

\medskip
Two particularly interesting choices  for the $w\in \bbr^m$ in Proposition 2.27 are
$$
w \ =\ (0,...,0,1,...,1) \qquad {\rm with\ } k \ {\rm ones,\ \ \  and}
\eqno{(2.32)}
$$
$$
w \ =\ (\d,...,\d,\d+1) \qquad {\rm with\ } \d>0.
\eqno{(2.33)}
$$
The first choice yields

\Prop{2.29. ($k$-Fold Sums)}  {\sl  Suppose $p$ is $a$-\hy of degree $m$ on $V$. For each 
$1\leq k\leq m$
$$
q(x) \ =\  {\prod_{|I|=k}}' \left( \l_{i_1}(x) +\cdots + \l_{i_k}(x)\right)
$$
defines an $a$-\hy polynomial of degree ${m\choose k}$ on $V$
whose eigenvalue functions are the $k$-fold sums}
$$
{1\over k}\left( \l_{i_1}(x) +\cdots + \l_{i_k}(x)\right)\qquad i_1 < \cdots < i_k.
$$

The second choice yields

\Prop{2.30. (``$\d$-Uniformly Elliptic'')}  {\sl  Suppose $p$ is $a$-\hy of degree $m$ on $V$. For each 
$\d>0$
$$
q_\d(x)\ =\ \prod_{k=1}^m  \left ( \l_k(x) + \d \tr \l(x)\right)
$$
defines an $a$-\hy polynomial of degree $m$ on $V$ whose eigenvalue functions are
${1\over m\d+1}$ times
$$
\l_k(x) + \d  \tr \l(x)  \qquad k=1,...,m.
$$
The G\aa rding cones $\G_\d$ for $q_\d$ form a   conically-fundamental
neighborhood system of the G\aa rding cone $\overline \G$ for $p$.
}

  
\vskip .2in 
\centerline{\headfont The Edge and the Polar}

\medskip

To begin we review some standard results for any  non-empty open convex cone $\G\ss V$ with vertex at the origin.  Set
$$
\G^+\ =\ \overline\G\and \G_+\ =\ {\overline\G}^0
\eqno{(2.34)}
$$
where the {\bf polar} $A^0$ of a closed convex cone $A$ with vertex 
at the origin is defined by
$$
A^0\ =\ \{y\in V : (x,y)\geq 0 \ \forall \, x\in A\}.
$$
Recall the Bipolar Theorem: $(A^{0})^0=A$. Note that $\span\G^+ = V$ since $\G$ is open and 
non-empty.

The {\bf edge} of $\G^+$ (or $\G$) can be defined as
$$
E\ =\ \{x\in V : x+\G^+ =\G^+\}
\eqno{(2.35)}
$$
One sees easily that $E$ is a convex cone contained in $\G^+$, and since 
$x\in E \iff -x\in E$ ($x+\G^+ =\G^+ \iff  \G^+ =-x+\G^+$), the edge $E$ must be a vector subspace of $V$.
Moreover,
$$
E\ {\rm contains \ any\ linear\ subspace\ of\ } \G^+.
\eqno{(2.36)}
$$
The edge can also be defined by
$$
E\ =\ \G^+ \cap (-\G^+).
\eqno{(2.37)}
$$
To see this recall that with $E$ defined by (2.35), $E\ss \G^+ \cap (-\G^+)$ was noted above.
If $x\in \G^+ \cap (-\G^+)$, then the line through $x$ is a subset of $\G^+$.  
By (2.36) this line is a subset of the edge $E$ defined by (2.35).

\medskip

Again for a  general open convex cone $\G$ we have:
$$
E\ {\rm and}\ \span \G_+ \ {\rm are\ polar\ cones}.
\eqno{(2.38)}
$$
{\bf Proof.} Note that $E\ss\G^+ \ \Rightarrow\ \G_+\ss E^0$ and hence $\span \G_+\ss E^0$.
Also we have $\G_+ \ss\span \G_+  \ \Rightarrow\ \ (\span \G_+)^0 \ss \G^+$. 
The polar of a vector space is a vector space, hence by (2.36) we have $(\span \G_+)^0\ss E$.
 Now apply the Bipolar Theorem.\qed

\medskip

If $E=\{0\}$, then $\G^+$ is said to be {\bf complete}.  This is equivalent to
 $V=\span \G_+$ or $\G_+ = \overline {\Int\G_+}$ by (2.38).
 \medskip
 
 Now we assume that $\G$ is the \Ga\  cone of an $a$-hyperbolic polynomial $p$.
 Then the edge $E$  of $\G(p)$  can be described in  several ways.

\Theorem{2.31} {\sl
Suppose  $p$ is an $a$-hyperbolic polynomial with \Ga\  cone $\G$.\smallskip

(a)\ The edge $E$ of $\G$ equals the null set of $p$, i.e., 
$$
E\ =\ \{x\in V :  \l_1(x) = \cdots = \l_m(x)=0\}.
\eqno{(2.39)}
$$

(b)\ The edge $E$ of $\G$ equals the {\sl linearity}  of $p$, i.e., }
$$
E\ =\ \{x\in V : p(b+tx) = p(b) \ \ \forall\, t\in \bbr, b\in V\}.
\eqno{(2.40)}
$$

(c) \centerline{$x\in E\quad\iff\quad p'_x\ \equiv\  0.$  \qquad\qquad\qquad\qquad}

\bigskip
\noindent
{\bf Proof of (a).} By 1(b) the negative of the \Ga\ cone, $-\G^+$,  is the set of points $x\in V$ whose $a$-eigenvalues are all
$\leq0$.  Thus $E=\G^+\cap (- \G^+)$ is the null set of $p$.
\smallskip

\noindent
{\bf Proof of (b) and (c).} Apply Property (3), and then for (b) note that
$ p(b+tx)-p(b)$ vanishes for all $b\in\G$ if and only if it vanishes for all $b\in V$.
\qed\medskip

This has several important consequences.
First, taking derivatives does not change the edge.

\Cor{2.32} {\sl
Suppose  $p$ is  $a$-hyperbolic of degree $m\geq 2$.  Then}
$$
{\rm Edge}(p_a')\ =\ {\rm Edge}(p).
$$
\pf
By (2.23) ${\rm Edge}(p_a') \supset{\rm Edge}(p)$.  If $\l_1'(x) = \cdots = \l_{m-1}'(x)=0$,
then all $\l_k(x)$ vanish by the sharp application (2.23)$'$ of Rolle's Theorem given in the proof 
of Proposition 2.21.  \qed
\medskip

Second, the edge of the restriction is the restriction of the edge.

\Cor{2.33} {\sl
Suppose  $p$ is  $a$-hyperbolic polynomial on $V$, and that $W$ is  a subspace of $V$.
If $\G\cap W \neq \emptyset$, then}
$$
\Edge\left(p\bigr|_W\right)  \ =\ \Edge(p)\cap W.
$$

\vskip .5in


\centerline{\headfont \ 3.\   $\G$-Monotone Sets}
\medskip

The Dirichlet problem will be studied (and solved) for equations determined
by subsets  $F$    of $V = \Symn$ which are $\G$-monotone.  In this section we discuss these subsets.
Recall that by Corollary 2.17  they  include all  branches of $\{p=0\}$.
Throughout this section we assume that $p$ is a given $a$-hyperbolic polynomial,
$\G$ is the associated G\aa rding cone, and $\l_1(x),...,\l_m(x)$ are the $a$-eigenvalue functions.

  
\vskip .2in 
\centerline{\headfont The Structure Theorem}

\smallskip

Recall (4) that
$$
{\rm trace}(x)\ =\ \l_1(x) +\cdots +\l_m(x).
\eqno{(3.1)}
$$
(i.e.,  the elementary symmetric function $\s_1(x)$)  is a linear function of $x$.  Define the hyperplane
$$
V_0\ \equiv\ \{x\in V : {\rm trace} (x)=0\}.
$$
We adopt the notation in Remark 2.14 with $W=V_0$.

\Def{3.1}  A function $f:V_0\to \bbr$ is {\bf  $\| \cdot\|^\pm$-Lipschitz} if
$$
-\|y\|^- \ \leq \ f(x+y)-f(x)\ \leq\ \|y\|^+ \fa x,y\in V_0.
\eqno{(3.2)}
$$

\Theorem{3.2. (The Structure Theorem)}  {\sl  Suppose $F$ is a closed subset of $V$ with $\emptyset \neq F\neq V$.
If $F$ is $\G$-monotone, then there exists a $\| \cdot\|^\pm$-Lipschitz function 
$f:V_0\to \bbr$ such that
$$
\partial F\ =\ \{ f(x) a +x : x\in V_0\}\ \ {\rm is\ the\ graph\ of\ } f,
\eqno{(3.3)}
$$
and
$$
 F\ =\ \{ t a +x : t\geq f(x) \ {\rm and\ } x\in V_0\}\ \ {\rm is\ the\ upper\ graph\ of\ } f.
\eqno{(3.4)}
$$
Conversely, each such $f$ defines a $\G$-monotone set $F$ via (3.4).}

\Cor{3.3}  {\sl Suppose $F$ is a closed subset of $V$ which is $\G$-monotone.
Then $\Int F$  is a connected open set with $F=\overline{\Int F}$.
Moreover, $F$ is uniquely determined by its boundary since $F=\partial F+\G$.}

\Ex{3.4. (Positive-Orphant-Monotone)}  
Take $V\equiv \bbr^m$ and $p(\l) \equiv \l_1\cdots\l_m$ where 
$\l=(\l_1,...,\l_m)\in\bbr^m$. Set $e=(1,...,1)$.
Then $p$ is $e$-hyperbolic with eigenvalues $\l_1,...,\l_m$ at $\l$.
The closure $\overline\G$  of the G\aa rding cone is the closed positive
orphant $\Or$.  The traceless hyperplane  $V_0$ is  the orthogonal
complement of $e$ while $\|\l\|^+ = -\l_{\rm min}$ and $\|\l\|^- = \l_{\rm max}$.
If $\E$ is a closed subset of $\bbr^m$ with $\emptyset\neq \E\neq \bbr^m$,
and $\E$ is $\Or$-monotone, then there exists $f:V_0\to \bbr$ which is 
$\| \cdot\|^\pm$-Lipschitz  continuous such that $\partial \E$ is the graph
of $f$ over $V_0$ and $\E$ is the upper graph of $f$ over $V_0$.
\medskip

The next example is the key example for studying the Dirichlet problem.

\Ex{3.5. (The Structure of $\cp$-Monotone Sets)}  
Take $V=\Symn$ and $p(A)=\det A$.  Then $p$ is $I$-hyperbolic with eigenvalues
exactly the standard eigenvalues of $A$.  The closed G\aa rding cone $\overline{\G} = \cp
\equiv \{A\geq0\}$.  Consider the hyperplane 
$$
\Sym_0(\rn) \equiv \{A\in \Symn: \tr A=0\}.
$$
If $F$ is a non-trivial closed subset of $\Symn$ 
($\emptyset\neq F \neq \Symn$) which is $\cp$-monotone, then there exists a function
$f:\Sym_0(\rn)\to \bbr$, which is $\|\cdot\|^\pm$-Lipschitz, such that 
$$
F\ =\ \{tI+A : A\in \Sym_0(\rn) \ {\rm and\ } t  \geq f(A)\}\ \ {\rm is\ the \ upper \ graph\  of
 \ } f.
 \eqno{(3.5)}
$$
Moreover, the properties of $F$ enumerated in Corollary 3.3 hold.

\medskip
\noindent
{\bf Proof of Theorem 3.2.}  The monotonicity hypothesis $F+\G\ss F$ implies that the 
set $\ell_x = \{t\in \bbr : ta+x\in F\}$ is either empty, all of $\bbr$, or of the form $[c,\infty)$
with $c\in\bbr$.  If $\ell_x=\emptyset$, then $F=\emptyset$ since the existence of $y\in F$
implies that $x+ta = y+(x-y+ta) \in F$ if $t$ is chosen large enough so that 
$x-y+ta \in \G$.  If $\ell_x=\bbr$, then $F=V$ since for each $y\in V$, if $t<<0$, then
$y-x-ta \in \G$ which implies that $y = x+ta +(y-x-ta) \in F+\G \ss F$.

Since $V$ is neither empty nor all of $V$, there is a well-defined function
$f:V_0\to\bbr$ such that $\ell_x = [f(x), \infty)$.
Note that $\partial F$ is the graph of $f$ over $V_0$.
This proves (3.3) and (3.4).

We now claim that if $F$ is  the upper graph of $f$ as in (3.4), then
\medskip
\centerline
{$F$ is $\G$-monotone  \quad $\iff$  \quad $F$  is  $\|\cdot\|^\pm$-Lipschitz.}
\medskip
To see this suppose $\bar x = f(x) a + x$, $x\in V_0$, is a point on the graph of $f$.
Then the cone $\bar x+\G$ lies above the graph of $F$ if and only if 
$$
f(x+y) \ \leq \ f(x) + \|y\|^+ \fa y\in V_0
$$
by Remark 2.14.  This inequality for all $y\in V_0$ is equivalent to
$$
-  \|y\|^- \ \leq \ f(x+y) - f(x)   \fa x, y\in V_0.
$$
\qed

\vskip .2in
 

\centerline{\headfont  A  Universal Construction of  $\G$-Monotone Subsets}
\smallskip

Let $\E \ss\bbr^m$  be a closed subset which is symmetric, i.e., invariant under the permutation of
coordinates on $\bbr^m$.
Then $\E$ universally determines a subset $F_\E = F_{\E}(p) \ss V$ 
for each $a$-hyperbolic polynomial $p$ of degree $m$ on $V$ by setting
$$
F_\E\ =\ \{x\in V : \l(x)\in \E\}\ =\ \l^{-1}(\E)
\eqno{(3.6)}
$$
where $\l(x)$ is the eigenvalue map (defined mod $\pi_m$).
If \E\ is a cone, then $F_\E$\ is a cone.
Note that the  G\aa rding cone $\overline \G$ for $p$ is universally determined by the positive orphant
$\Or$.  The next three theorems will be incorporated into one of our main results
(Theorem 5.19) on the Dirichlet problem. The most basic for this application is the following.

\Theorem {3.6}  {\sl  Suppose $\E$ is a closed symmetric subset of $\bbr^m$.
\smallskip
\centerline{If  $\E$  is $\Or$-monotone, then $F_\E$ is $\G$-monotone.}
}

\pf
Suppose $x\in F_\E$, i.e., $\l_a(x)\in \E$.
Pick $b\in\G$.  Using the real analytic arrangement Theorem
2.9  states that the vector $\l^+$ with coordinates
$\l_k^+ \equiv \l^k_a(x+b)-\l^k_a(x)$ satisfies
$$
\l^+\ \in\ \Or.
$$
Hence $\l_a(x+b) =  \l_a(x) +\l^+ \in \E +\Or \ss\E$ proving that $x+b\in F_\E$.\qed

\medskip

Given an $a$-hyperbolic polynomial $p$ of degree $m$ with eigenvalue map
$\l(x)$ and \Ga \  cone $\G$, consider the conical neighborhoods $\overline{\G_\d}$
of $\overline \G$ defined by 
$$
\l_k(x)  + \d \tr \l(x)  \ \geq\ 0\qquad k=1,...,m
$$
as in Proposition 2.30.  Let $\left( \Or \right)_\d \ss \bbr^m$ denote the conical neighborhood of the 
positive orphant $\Or$   defined by 
$$
\l_k  + \d \tr \l  \ \geq\ 0\qquad k=1,...,m
$$

\Theorem{3.7}  {\sl Suppose that  $\E$ is a closed symmetric subset of $\bbr^m$
and $\d>0$. 
\medskip
\centerline
{ If  $\E$ is $\left( \Or \right)_\d$-monotone, then $F_\E$ is $\overline{\G_\d}$-monotone.
}
}

\pf
Consider the $a$-hyperbolic polynomial $q_\d$ defined in Proposition 2.30.  It has 
eigenvalue map (up to scale) 
$
\L(x) \ =\ \l(x) + (\d \tr \l(x)) e,
$
and \Ga\ cone $\overline{\G_\d} = \L^{-1}(\Or)$.  In order to apply Theorem 3.6 to $q_\d$ consider the
linear map $L:\bbr^m \to \bbr^m$ defined by $L(\l) = \l +(\d\tr\l)e$.  Note that 
$L^{-1} (\L) = \L - \left( {\d\over 1+m\d} \tr \L   \right)e$.  It suffices to show that:
$$
F_\E\ =\ \L^{-1}(L\E), \ \ {\rm and\ that}
$$
$$
\E\ {\rm \ is\ \ } \left( \Or \right)_\d\ {\rm monotone} \qquad\iff\qquad L\E\ \ {\rm is\ \ } \bbr^m_+ \  {\rm monotone}
\eqno{(3.7)}
$$
to conclude by Theorem 3.6 that $F_\E = \L^{-1}(LE)$ is $\overline \G_\d= \L^{-1}(\Or)$-monotone.
First, $\L(x) = L\l(x)$ implies $F_\E =   \l^{-1}(\E) =  \L^{-1}(L\E)$. Second, $ \left( \Or \right)_\d = L^{-1} \Or$ implies
(3.7). \qed
\medskip

Finally we prove a  recent convexity result of Bauschke, G\'uler, Lewis and Sendov [BGLS]
stated first with our applications in mind.

\Theorem{3.8}  {\sl Suppose that \E\ is a symmetric subset of $\bbr^m$.
\medskip
\centerline
{
If \E\ is convex, then $F_\E$ is convex.
}
}

\pf
Suppose $x\notin F_\E$, i.e., $\lu(x) \notin \E$. It will suffice to exhibit a convex set $A$
which contains $F_\E$ but not $x$.

Note that by Corollary 2.28, for each $\wu\in \bbr^m_{\ua}$, the set 
$$
A\ \equiv \ \{y\in V : \wu(\lu(y)) <1\}
$$
is convex.  By the Hahn-Banach Theorem there exists a $w\in (\bbr^m)^*$ such that
$$
\sup_{\E} w \ <\ 1 \ = \ w\left( \lu(x)\right).
$$
Since $\sup_{\E} \wu = \sup_{\E} w$ and $w\left(\lu(x)\right) \leq \wu\left(\lu(x)\right)$,
$$
\qquad\qquad\qquad\qquad\qquad\qquad
\sup_{\E} \wu \ <\ 1\ \leq \wu\left(\lu(x)\right). \qquad\qquad\qquad\qquad\qquad\qquad \mathqed
$$

\Remark{3.9}  Adopt the notation in Corollary 2.28 and set $u\equiv \lu(sx+(1-s)y), v= s\lu(x)+(1-s)\lu(y)$.
Then for $\wu\in \bbr^m_{\ua}$ the Corollary concludes that $\wu(u) \leq \wu(v)$. Now
$$
u\ \in\ \E \ \equiv\ {\rm the\ convex\ hull\ of\ the\ orbit \ } \pi_m v.
\eqno{(3.8)}
$$
Otherwise, as in the proof of Theorem 3.8, there would exist $\wu\in  \bbr^m_{\ua}$ with
$$
\sup_{\E} \wu \ <\ 1\ \leq\ \wu(u)
$$
which contradicts $\wu(u) \leq \wu(v)$.   

\Theorem{3.10} {\sl 
Suppose $f$ is a symmetric function on $\bbr^m$.}
$$
\eqalign
{
&\ \ {\sl If \ } f(\l)\ {\sl is\ convex, \ then\ } f(\l(x)) \ {\sl is\ convex}
}
\eqno{(3.9)}
$$
\pf Note that by (3.8) there exist $v_1,...,v_{N}$  
in the orbit $\pi_m v$ and $0\leq t_j\leq 1$, $j=1,...,N$ with $\sum_j t_j=1$
such that $u=\sum_j t_j v_j$.  By the convexity and symmetry of $f$,
$$
f(u)\ =\ f\left( \sum_j t_j v_j\right) \ \leq\ \sum_j t_j f( v_j) \ = \ \sum_j t_j f(v) \ =\ f(v).  
\qquad\qquad \vrule width5pt height5pt depth0pt
$$

\vskip .3in
 

\centerline{\headfont \ 4.\  The Dirichlet Problem}
\medskip

In this section we summarize without proof the existence and uniqueness theorem in  [HL$_1$]  
for the Dirichlet problem.  It can be read independently of Sections 2 and 3.
The equations considered are on euclidean space.  They are purely second order
and have constant coefficients; but are fully nonlinear and degenerate elliptic.

\vskip .2in

\centerline{\headfont Subequations and Positivity}
\medskip
Each subset $F\ss \Symn$ determines a class of $C^2$-functions $u$ on an open subset
$X\ss \rn$ by the requirement that its hessian  (second derivative)  $D^2u(X)$ lie in $F$ at each $x\in X$.
This class in $C^2(X)$ can be extended to a class in $\USC(X)$, the upper semi-continuous
$[-\infty, \infty)$-valued functions on $X$, provided that $F$ satisfies monotonicity with
respect to the set
$$
\cp\ \equiv\ \{A\in\Symn : A\geq0\}
\eqno{(4.1)}
$$
Two justifications for imposing this condition are given in Remark 4.3 and Remark 4.10.

\Def{4.1}  A closed subset $F\ss\Symn$ with $\emptyset\neq F\neq \Symn$ will be 
called a   {\bf subequation} or {\bf Dirichlet set} if $F$ satisfies the {\sl positivity condition}
$$
F + \cp \ \ss\ F.
\eqno{(4.2)}
$$

The Structure Theorem 3.2 for $\cp$-monotone sets, as explained  in Example 3.5,
is now the ``Structure Theorem for Subequations''. In particular, Corollary 3.3 states that
$\Int F$ is  connected,  $F=\overline{\Int F}$, and $F$ is completely determined
by its boundary. In fact, it is the upper graph of a Lipschitz function on $\Sym_0(\rn)
=\{\tr A=0\}$.

\Def{4.2}  Suppose that $F\ss\Symn$ is a subequation
and  $u\in C^2(X)$ (where $X$ is an open subset of $\rn$).
If $D^2u(x)\in F$ for all $x\in X$, then $u$ is called {\bf $F$-subharmonic on} $X$.
If $D^2u(x)\in \Int F$ for all $x\in X$, then $u$ is called {\bf strictly $F$-subharmonic on} $X$.
If $D^2u(x)\in \partial F$ for all $x\in X$, then $u$ is called {\bf $F$-harmonic on} $X$.

\Remark{4.3. (Motivation)}
This definition of $F$-subharmonicity will be extended to non-differentiable functions. Perhaps the 
most fundamental property is that the maximum $\max\{u,v\}$ of two $F$-subharmonic functions
$u$ and $v$ is again $F$-subharmonic.  Computing the distributional hessian of
 $\max\{u,v\}$ in the case where $u$ and $v$ are $C^2$ provides strong support for assuming
 positivity if one wants $\max\{u,v\}$ to be $F$-subharmonic.  (See Remark 3.3 in  [HL$_1$]).

\bigskip

\centerline{\headfont Dirichlet Duality}
\medskip

Each subequation $F$ has a natural ``dual'' subequation.

\Def{4.4}  The {\bf Dirichlet dual} of a subequation $F\ss\Symn$ is the set
$$
\ft\ =\ \sim(-\Int F) \ =\ -(\sim \Int F).
\eqno{(4.3)}
$$

\Prop{4.5}  {\sl  The Dirichlet dual $\ft$ of a subequation $F$ is also a subequation.
Furthermore, duality holds, that is}
$$
\widetilde{\ft}\ =\ F.
\eqno{(4.4)}
$$
\pf
Assertion (4.4) follows straightforwardly from the fact that $F=\overline{\Int F}$.
To prove that $\ft$ satisfies (4.2) a lemma is required.

\Lemma{4.6}  {\sl Suppose $F=\overline{\Int F}$.  Then }
$$
\wt{F+A} \ =\ \ft-A \quad  {\rm for\ each\ } A\in\Symn.
\eqno{(4.5)}
$$
\pf
Note that $B\in \wt{F+A} \ \iff\ -B \notin \Int(F+A) = \Int F + A \ \iff\ -( B+A) \notin \Int F
\ \iff\  B+A \in \ft$.\qed

\medskip\noindent
{\bf Proof of Proposition 4.5.}
Suppose $P\in\cp$.  Then  $F+P\ss F$, or equivalently $F\ss F-P$. 
Taking duals reverses this inclusion and by Lemma 4.6, 
$\wt{F-P} =\ft+P$.  Thus $\ft+P\ss\ft$ as desired.\qed

\Prop{4.7}  {\sl Suppose $F$ is a subequation and $u\in C^2(X)$.
Then 
\smallskip
\centerline{
$u$ is $F$-harmonic \ \ $\iff$\ \  $u$ is $F$-subharmonic
and $-u$ is $\ft$-subharmonic.}}

\pf 
Note that $\partial F = F\cap (\sim\Int F) = F\cap (-\ft)$.\qed

\bigskip

\centerline{\headfont F-Subharmonic and F-Harmonic Functions}
\medskip

It is necessary and quite useful to extend the  definition of $F$-subharmonic 
to non-differentiable functions $u$.
Let $\USC(X)$ denote the set of $[-\infty, \infty)$-valued, upper semicontinuous 
functions on $X$.

\Def{4.8} A function $u\in \USC(X)$ is said to  be  $F$-{\bf subharmonic} 
if for each $x\in X$ and each function $\vf$ which is $C^2$ near $x$, one has that
$$
\left.
\cases
{ u-\vf  \ &$\leq$ \ 0 \    \quad {\rm near}\ $x_0$  and \cr 
  \ &= \ 0\ \qquad {\rm at}\ $x_0$  
 } 
 \right\} 
 \qquad \Rightarrow\qquad
 D^2_x  \vf \ \in \ F.
\eqno{(4.6)}
$$

Note that if $u\in C^2(X)$, then
$$
u\in F(X)\quad\Rightarrow\quad  D^2_x u \in F \qquad \forall \ x\in X
$$
since the test function $\vf$ may be chosen equal to $u$ in (4.10).  The converse   is not
  true for general subsets $F$.  However, we have the following.

\Prop{4.9} {\sl Suppose  $F$ satisfies the  Positivity Condition (4.2) and $u\in C^2(X)$.
Then}
$$
D^2_x u \ \in\ F \fa x\in X \qquad\Rightarrow\qquad u\ \in\ F(X).
$$

\pf
Assume $D^2_{x_0} u \in F$ and  and $\vf$ is a $C^2$-function such that  
$$
\left.
\cases
{ u-\vf  \ &$\leq$ \ 0 \    \quad {\rm near}\ $x_0$   \cr 
  \ &= \ 0\ \qquad {\rm at}\ $x_0$  
 } 
 \right\} 
$$
Since $(\vf - u)(x_0) = 0$,
$D{x_0} (\vf-u)_ = 0$, and $\vf-u \geq 0$ near $x_0$,  we have $D^2_{x_0}(\vf-u) \in \cp$.
 Now the Positivity Condition implies that  
 $D^2_{x_0}\vf \in  D^2_{x_0} u +\cp \ss F$.  This proves
that $u\in F(X)$.   \qed  \medskip

\Remark {4.10} Because of Proposition 4.9 {\sl we must assume that
$F$ satisfies the Positivity Condition}. Otherwise the definition of  $F$-subharmonicity would  not
extend the natural one for smooth functions (cf. Remark 4.3). 
The positivity condition is rarely used in proofs.  This is because without it $F(X)$ is empty and the results are trivial. For example, the positivity condition is not required in the following theorem. ($F$ need only be closed.)
\medskip

It is remarkable, at this level of generality,   that $F$-subharmonic
 functions  share many of the important properties
of classical subharmonic functions.   

\Theorem{4.11. Elementary  Properties of  F-Subharmonic Functions}
{\sl
Let $F$ be an arbitrary  closed subset of $J^2(X)$.
\medskip

\item{}  (Maximum Property)  If $u,v \in F(X)$, then $w=\max\{u,v\}\in F(X)$.

\medskip

\item{}     (Coherence Property) If $u \in F(X)$ is twice differentiable at $x\in X$, then $\jtx u\in F_x$.

\medskip

\item{}  (Decreasing Sequence Property)  If $\{ u_j \}$ is a 
decreasing ($u_j\geq u_{j+1}$) sequence of \ \ functions with all $u_j \in F(X)$,
then the limit $u=\lim_{j\to\infty}u_j \in F(X)$.

\medskip

\item{}  (Uniform Limit Property) Suppose  $\{ u_j \} \ss F(X)$ is a 
sequence which converges to $u$  uniformly on compact subsets to $X$, then $u \in F(X)$.

\medskip

\item{}  (Families Locally Bounded Above)  Suppose $\cf\subset F(X)$ is a family of 
functions which are locally uniformly bounded above.  Then the upper semicontinuous
regularization $v^*$ of the upper envelope 
$$
v(x)\ =\ \sup_{f\in \cf} f(x)
$$
belongs to $F(X)$.

}

\pf
 See Appendix B in [HL$_2$].

\vskip .5in

\centerline{\headfont Boundary Convexity}
\medskip

Associated to each subequation $F$ is an open cone $\Fa$ which governs the geometry 
of those domains for which the Dirichlet problem is always solvable. 
If $F$ is a cone with vertex at the origin, then $\Fa' = \Int F$.  Otherwise,  $\Fa' $ is constructed
as in [HL$_1$, pp. 415-416]. Here is a summary.

Fix a vertex $B\in\Symn$ and consider the ray sets:
$$
\overrightarrow{F_B} \ =\ \{A\in\Symn : B+tA \in F \ \forall \, t \geq \ {\rm some}\ t_0\}
\eqno{(4.7)}
$$
$$
\overrightarrow{F_B}' \ =\ \{A\in\Symn : B+tA \in \Int F \ \forall \, t \geq \ {\rm some}\ t_0\}
\eqno{(4.8)}
$$
Examples show that $\overrightarrow{F_B}$ may not be closed, 
$\overrightarrow{F_B}'$ may not be open,
and $\overrightarrow{F_B}, \overrightarrow{F_B}'$ may depend on $B$.  However,
$$
\Fa \ \equiv \ {\rm closure} \overrightarrow{F_B}
\and 
\Fa' \ \equiv \  \Int \overrightarrow{F_B}'
\eqno{(4.9)}
$$
are independent of the choice of vertex $B$.  In addition
$$
\Fa \  {\rm is\  a\  subequation\ which\ is \  a\  cone\  with\  vertex\  at \ the\ origin.}
\eqno{(4.10)}
$$
Moreover, by Lemma 5.8
 and Elementary Property (3),
$$
\Fa' \ = \ \Int \overrightarrow{F}
\and 
\Fa \ \equiv \  {\rm closure} \overrightarrow{F}'
\eqno{(4.11)}
$$

\Def{4.12} We call $\Fa'$ the {\bf asymptotic interior} of $F$, and $\Fa$ the 
{\bf asymptotic subequation} of $F$. (In [HL$_1$], $\Fa$ is called the ``asymptotic ray set associated to $F$''.)
\medskip

Note that (4.10) and (4.11) imply that 
$$
\Int \cp \ \ss\ \Fa'\ =\ \Int \Fa.
\eqno{(4.12)}
$$

Suppose now that $\O\ss\ss \rn$ is a domain with smooth boundary $\bo$.
Let $\rho$ be  a {\sl local defining function} for $\bo$ at a point $x$, that is, a smooth function
 with $|\nabla \rho|>0$, defined on a neighborhood $U$ of $x$, such that
$U\cap\O = \{\rho<0\}$

\Def{4.13}  The boundary $\bo$ is {\bf strictly $\Fa$-convex at $x$  if} 
$$
\Hess_x \rho \bigr|_{T_x\bo} \ =\ B\bigr|_{T_x\bo} \qquad{\rm for\ some\ } B\in \Fa' \equiv \Int \Fa. 
\eqno{(4.13)}
$$
This  is independent of the choice of  $\rho$ 
and is equivalent to the condition that:
$$
\Hess_x \rho+ tP_n\  \in\ \Fa' \equiv \Int \Fa  \qquad\forall\, t\ \geq \ {\rm some\ }t_0
\eqno{(4.14)}
$$
where $P_n$ denotes orthogonal projection onto the normal line to $\bo$ at $x$.
\medskip

Note that each classically strictly convex boundary is strictly $\Fa$-convex by 
(4.12).  In general
$\Fa$-convexity can be expressed purely in terms of the geometry of the boundary.
Let $II$ denote the second fundamental form of $\bo$ with respect to the interior
normal at $x$.  Then $\bo$ is strictly $\Fa$-convex at $x$ if and only if either of
the following equivalent conditions holds.
$$
\eqalign
{
&II \ =\ B\bigr|_{T_x\bo} \qquad{\rm for\ some\ } B\in \Fa  \cr
&II + tP_n\  \in\ \Fa' \equiv \Int \Fa  \qquad\forall\, t\ \geq \ {\rm some\ }t_0
}
$$
The concept of $\Fa$-convexity is important because it leads to the existence of
local barriers used in existence proofs. For domains $\O\ss\ss\rn$ we have the following global result.

\Theorem {4.14} {\sl Suppose $\bo$ is strictly $\Fa$-convex at each point.
Then there exists a global defining function $\rho\in C^\infty(\ob)$ and constants $\e_0>0, R_0>0$ such
that
\medskip
\centerline{
$C\left( \rho-\e\half |x|^2\right)$ is strictly $F$-subharmonic on $\O$ for all $C\geq C_0$ and $0\leq\e\leq\e_0$.
}
}
\medskip

See  [HL$_1$, \S 5] for proofs of the above.


\vskip .2in

\centerline{\headfont Existence and Uniqueness}
\medskip

\Theorem {4.15}  {\sl  Suppose $F$ is a subequation in $\Symn$ 
and $\O$ is a bounded domain in $\rn$ with smooth boundary $\bo$.
If $\bo$ is both $\Fa$ and $\overrightarrow{\ft}$ strictly convex, then
{\rm the Dirichlet problem is uniquely solvable for all continuous boundary data}.
That is, for each $\vf\in C(\bo)$, these exists a unique $u\in C(\ob)$ satisfying:
\medskip

(1)\ \ $u$ is $F$-harmonic on $\O$, and
\medskip
(2)\ \ $u=\vf$ on $\bo$.
}
\medskip

Since $\cp\ss \ft$ and $\cp\ss \overrightarrow{\ft}$, if $\bo$ is strictly convex
(in the standard way), then
it is both $\Fa$ and $\overrightarrow{\ft}$ strictly convex. In particular,
the Dirichlet problem is uniquely solvable on all domains with strictly
convex boundaries, such as a ball.

This Theorem 4.15 was first proved in [HL$_1$].  The proof involved the use of
``subaffine'' functions (equivalently $\cpt$-subharmonic functions) and relied on a deep result of
Slodkowski on quasi-convex functions [S].  The following definition of $F$-subharmonic
functions was adopted there.

\Def{4.8$'$}  A function $u\in \USC(X)$ is said to be {\sl $F$-subharmonic} if for 
each $x\in X$ and each $C^2$-function $v$ which is $\ft$-subharmonic 
near $x$, the sum $u+v$ is subaffine near $x$,
\medskip

See Remark 4.9 in [HL$_1$] for a proof that Definitions 4.9 and 4.9$'$ are equivalent.

A second proof of Theorem 4.15 was given in [HL$_2$] using standard viscosity methods,  i.e., 
the Theorem on Sums and Definition 4.8 (see  [C], [CIL]).

\vfill\eject
 

\centerline{\headfont 5.   Subequations Determined by Hyperbolic Polynomials}
\medskip

G\aa rding's theory provides a unified approach to studying many basic subequations.
For this application to subequations we only consider   polynomials
on the vector space $V\equiv \Symn$ which are  hyperbolic in the direction $I$.

\vskip .2in

 \centerline{\headfont \DG Polynomials}

\Def{5.1}  A homogeneous real polynomial $M$ of degree $m$ on $\Symn$ is  
{\bf  $I$-hyperbolic}  if $M(I)>0$ and
 for all  $A\in \Symn$ the polynomial $M_A(s) = M(sI+A)$
 has $m$ real roots.
 \medskip

We adopt the terminology and notation from the previous discussion, with $\G$ denoting
the G\aa rding cone associated to the polynomial $M$.
The principal branch (or {\bf closed G\aa rding cone}) $F_1=\overline\G$
is defined by $\l_k(A)\geq 0, k=1,...,m$, and is a convex cone (Corollary 2.13) but
may not be a subequation.  We require that $\overline \G$ be a subequation.

\Def{5.2}  An $I$-hyperbolic polynomial  $M$  on $\Symn$
is said to be a {\bf \DG} polynomial if $\overline\G$ is a subequation, i.e., if positivity holds:

\medskip

\qquad
{(1)}  \ \ $\overline\G+\cp\ \ss\ \overline\G$.

\Remark{5.3} There are other useful ways of stating this positivity condition for $\overline \G$.
Since $\overline\G$ is a convex cone, (1) is equivalent to 

\medskip

\qquad
{(2)}  \ \ $\cp\ \ss\ \overline\G$, \quad that is,\ \ \   $A\geq 0 \ \ \Rightarrow\ \ \l_k(A) \geq0, \ \ k=1,...,m.$
\medskip
In other words, for $A\geq 0$,  all the roots of the polynomial  $t\mapsto M(tI+A)$
are $\leq 0$.  Equivalently,
\medskip

\qquad{(2)$'$}  \ \ $M(tI+A)\  \neq \ 0$ if $A\geq0$ and $t>0$.
\medskip
\noindent
Since the extreme rays in $\cp$ are generated by orthogonal projections $P_e$ 
with $e\in\rn$, it is enough to verify (2) for $A=P_e$.
\medskip
\qquad{(3)}  \ \ $\l_k(P_e)\ \geq\ 0$ for all unit vectors $e\in \rn$ and $k=1,...,m$.
\medskip
\qquad(3)$'$\ \ $M(tI+P_e)\ >\ 0$ for all $t>0$ and all unit vectors $e\in\rn$.

\medskip
\noindent
By definition the branches $F_k$ are subequations if  positivity holds. That is,
\medskip
\qquad(4)\ \ $F_k +\cp\ \ss\ F_k$ for all $k=1,...,m$.
\medskip
\noindent
By Corolloary 2.17  and the convexity of the principal branch $\overline\G = F_1$, condition (4) is equivalent
to condition (2).  Finally, by the Structure Theorem 3.2 the  condition that:
\medskip
\qquad(5)\ \ The  ordered eigenvalue function $\lu_k(x)$ is $\|\cdot\|^\pm$-Lipschitz ($k=1,...,m$),
\medskip
\noindent
is equivalent to (4).

\vskip .2in

 \centerline{\headfont Branches Considered as Subequations}
\smallskip

Some of G\aa rding's results from Section 2
can be summarized as follows.  The ordered eigenvalues determine $m$ branches
$F_k \equiv \{\lu_k\geq0\}$ of $\{M=0\}$.  Each $\partial F_k = \{\lu_k=0\}$ is contained in $\{M=0\}$,
and $F_k = \overline{\Int F_k}$ with $\Int F_k = \{\l_k>0\}$.  Each $F_k$ is a cone,
and the principal branch $\overline\G = F_1$ 
 is a convex cone, which is a monotonicity cone for each of the
other branches, i.e., 
$$
F_k+\overline\G\ \ss\ F_k\qquad k=1,...,m.
\eqno{(5.1)}
$$

\Theorem{5.4}  {\sl  Suppose that $M$ is a \DG polynomial of degree
$m$ on $\Symn$.  Then each branch $F_k$, $k=1,...,m$, is a subequation 
with dual subequation $\wt F_k = F_{m-k+1}$.
Moreover, each branch is a cone.}

\pf Property (1b), that $\lu_k(-A) = -\lu_{m-k+1}(A)$, for the ordered eigenfunctions, implies that 
$\wt F_k = F_{m-k+1}$. Property (1a), that $\l_k(tA) = t\l_k(A)$ for $t\geq0$, implies  $F_k$ is a cone.\qed

\medskip

Thus, as a special case of Theorem 4.15 we can solve the 
Dirichlet problem for each branch of the equation $\{M=0\}$.

\Theorem{5.5}  {\sl  Suppose that $M$ is a \DG \ polynomial on $\Symn$.
Let $i\equiv\max\{k, n-k+1\}$.  Then for each bounded domain $\O$ with smooth
strictly $\oa{F_i}$-convex boundary, the Dirichlet problem for $F_k$-harmonic
functions is uniquely solvable for all continuous boundary data.}

\Remark{5.6} (a) \ \ Fix $c\in\bbr$ and replace $F_k=\{\lu_k\geq0\}$ by 
 $F_k^c =\{\lu_k\geq c\}$.  Both Theorem 5.4 and Theorem 5.5 remain true except
 for the fact that $F^c_k$ is no longer a cone.  Note that $F^c_1$ remains a convex subequation.\smallskip
 
 (b) \ \ The principal branch $F_1=\overline\G$ can be perturbed to the convex cone subequation
 $(c>0)$
 $$
 \overline{\G_c} \ =\ \{ A\in \G : M(A)\geq c\}.
 $$
Convexity follows from Corollary 2.13 (d) or (e).

\vskip .2in


\centerline{\headfont  Examples of  \DG  Polynomials}
\medskip

 There are a number of geometrically interesting nonlinear equations that 
 arise from the construction in Theorem 5.4, and for which Theorem 5.5 solves the Dirichlet problem.
 
 \Ex{5.7. (The Real Determinant)}  Of course the fundamental example of a \DG
polynomial is the polynomial  on $\Symn$ given by
 $$
 M(A)\ =\ \det A
 \eqno{(5.2)}
 $$
 where the $M$-eigenvalues of $A$ are just the usual eigenvalues of $A$
 as a symmetric matrix.    The corresponding
 $F_1$-harmonic functions are convex solutions  $u$  of the homogeneous Monge-Amp\`ere
 equation  $\det \Hess\, u =0$ (in the viscosity sense (cf.  [HL$_1$], [HL$_2$])). The other branches
 $F_k$, $k=2,...,n$  lead to other branches of this equation.
 Theorem 5.5 says that the Dirichlet problem is uniquely solvable for $F_k$-harmonic
 functions on any domain which is strictly $\overrightarrow{F_i}$-convex
 where $i=\max\{k, n-k+1\}$.

 \Ex{5.8. (The Complex and Quaternionic Determinant)}  Consider $\bbc^n = (\bbr^{2n},J)$ where 
 $J:\bbr^{2n} \to \bbr^{2n}$ with $J^2=-I$ is the standard almost complex structure.
 Then any $A\in \Sym(\bbr^{2n})$ has a {\sl hermitian symmetric component}
 $$
 A_\bbc\ \equiv\ \half(A - JAJ)
 $$
 which commutes with $J$.  The   eigenvalues of $A_\bbc$,  considered
 as a real  $2n\times 2n$ symmetric matrix,  occur in pairs, $\l_1, \l_1, \l_2,\l_2,...,
 \l_n,\l_n$ with complex lines as eigenspaces, while $A_\bbc$,  considered
 as an $n\times n$  complex hermitian symmetric matrix,
  has $n$ real eigenvalues $\l_1(A_\bbc)=\l_1,...,\l_n(A_\bbc)=\l_n$.

 The complex determinant
 $$
 \det_\bbc(A_\bbc) \  = \ \l_1(A_\bbc)\cdots \l_n(A_\bbc)
 \eqno{(5.3)}
 $$
 is a \DG polynomial in $A\in \Sym(\bbr^{2n})$.
  \medskip
 
 Similarly, one can consider  quaternion $n$-space $\bbh^n = (\bbr^{4n}, I, J, K)$,
where $I,J,K$ satisfy the usual quaternion relations.  Then any $A\in \Sym(\bbr^{4n})$ has a 
{\sl  quaternionic  hermitian symmetric component}
 $$
 A_\bbh\ \equiv\ \smfrac 1 4 (A - IAI-JAJ - KAK)
$$
 which commutes with $I,J$ and $K$.  The  eigenvalues of $A$,  considered
 as a real  $4n\times 4n$ symmetric matrix,  occur in  multiples of four $\l_1,\l_1,\l_1,\l_1,\l_2,...$
  with quaternion lines as eigenspaces,
 while $A_\bbh$,  considered as a quaternionic hermitian symmetric $n\times n$  matrix can be diagonalized 
 under the action of ${\rm Sp}_n\cdot{\rm Sp}_1$ with $n$ real eigenvalues 
 $\l_1(A_\bbh)=\l_1,...,\l_n(A_\bbh)=\l_n$.
The  quaternion determinant
 $$
M(A)\ =\  \det_\bbh(A_\bbh) \  = \ \l_1(A_\bbh)\cdots \l_n(A_\bbh)
 \eqno{(5.4)}
 $$
is a \DG polynomial in $A\in \Sym(\bbr^{4n})$.

 \Ex{5.9. (Lagrangian Harmonicity)}  This is a non-classical case.
 Consider $\bbc^n = (\bbr^{2n},J)$ as above, and for $A\in \Sym(\bbr^{2n})$
 define the {\sl skew hermitian component}
 $$
 A_{\rm skew}\ \equiv\ \half(A +JAJ)
 $$
 which anti-commutes with $J$.  The eigenvalues occur in opposite pairs $\pm\mu_1,...,\pm\mu_n$
 where each $\pm \mu_k$-eigenspace generates a complex line.  Let $\tau= \half {\rm trace}_\bbr A$.
 
 \Prop{5.10}  {\sl The product
 $$
 M_{\rm LAG}(A) \ \equiv\ \prod_{ 2^n\ {\rm times}} \left( \tau \pm \mu_1\pm \cdots\pm\mu_n\right)
 \eqno{(5.5)}
 $$
 taken over all sequences $\pm\cdots\pm$, is a \DG polynomial on $\Sym(\bbr^{2n})$.}

 More generally, for any $1\leq p\leq n$, the product
 $$
 M(A)\ =\  {\prod_{|I|=p}}' \left(  \tau \pm \mu_{i_1} \pm \cdots\pm  \mu_{i_p}  \right)
  \eqno{(5.6)}
 $$
 is a \DG polynomial.  (See [HL$_3$] for the proof.)
 \medskip
 
 We shall refer to these four basic examples as the {\sl real, complex, quaternionic, and 
 Lagrangian/Isotropic Monge-Amp\`ere polynomials}.  Now we describe three methods of constructing
 a new \DG polynomial from a given one.
 
 The first method decreases the degree of the \DG polynomial.
 
 \medskip
 \noindent
 {\bf Method 5.11. (Derivatives -- The Elementary Symmetric Functions).}
 Given a \DG polynomial $M$ with \Ga\ cone $\G$, for each $A_0\in\G$, the $k$th derivative
 $M^{(k)}$ in the direction $A_0$ is also a \DG polynomial.  This is because 
 $\G_{M^{(k+1)}} \supset \G_{M^{(k)}} \supset   \cdots  \supset \G_{M}$,
 by Proposition 2.21, and so $\cp\ss\overline\G_{M}$
 implies   $\cp\ss\overline\G_{M^{(k)}}$.
 By (2.26) this construction is equivalent to taking the elementary symmetric functions
 of the $A_0$-eigenvalues (as functions on $\Symn$).
 
 This method can be applied, of course, to the real, complex, quaternionic, and 
 Lagrangian/Isotropic Monge-Amp\`ere polynomials.
 
 \medskip
 \noindent
 {\bf Note.}  Different choices of $A\in \G$ produce different \DG polynomials.
 \medskip
 
 The second method increases the degree of the \DGP.

 \medskip
 \noindent
 {\bf Method 5.12. ($k$-Fold Sums of Eigenvalues).}
 Suppose we are given a \DGP\ $M$ on $\Symn$ with \Ga\  cone $\G$.
 Fix $A_0\in \G$ and consider the $A_0$-eigenvalue functions $\l_1(A),...,\l_m(A)$. Then
 $$
 M_k(A)\ =\    {\prod_{|I|=k} }'    \left  (\l_{i_1}(A)+\cdots +\l_{i_k}(A)  \right)\qquad A\in \Symn
  \eqno{(5.7)}
 $$
 defines a \DGP\ of degree ${m\choose k}$ on $\Symn$ with \Ga\ cone $\G_k\supset \G$
 and $A_0$-eigenvalues functions
 $$
\L_I(A)  \ =\ {1\over k}\left(  (\l_{i_1}(A)+\cdots +\l_{i_k}(A)\right).
  \eqno{(5.8)}
 $$
This follows from Proposition 2.29 and Theorem 3.6.

\Note{5.13} 
\smallskip

 (a) \ It is easy to see that $\G= \G_1\ss\G_2 \ss\ \cdots\ \ss\G_n$.
\smallskip

(b) \ Any vector $w\in\Int \Or$ can replace $e=(1,...,1)$.  We would then have:
$$
 M_k(A)\ =\    {\prod_{|I|=k} }'    \left  ( w_{i_1}\l_{i_1}(A)+\cdots + w_{i_k}\l_{i_k}(A)  \right)\qquad A\in \Symn
 $$
 
 The third method leaves the degree of the \DGP\  fixed.

 \medskip
 \noindent
 {\bf Method 5.14. ($\d$-Uniformity).}
Define for $\d>0$
$$
\cp_\d\ =\ \{A\in\Symn : A+(\d \tr A)\cdot I\geq 0\}.
\eqno{(5.9)}
$$
Each $\cp_\d$ is a convex cone subequation whose interior contains $\cp$.
the family is {\sl conically-fundamental} for $\cp$ in the following sense. 
$$
\eqalign
{
&{\rm If\ } F\ {\rm is\ any\ convex\ cone\ subequation\ with\ }  \cp-\{0\} \ss\Int F,   \cr
&\ \ \ \ {\rm then\ there\ exists\ } \d>0\ {\rm such\ that\ }\cp\ss\cp_\d\ss F
}
\eqno{(5.10)}
$$

\Def{5.15}  A subequation $F\ss\Symn$ is {\bf uniformly elliptic} if for some 
$\d>0$, $F$ is $\cp_\d$-monotone, i.e., 
$$
F+\cp_\d\ \ss\ F.
\eqno{(5.11)}
$$
Using (5.10)  it is easy to show that our definition of uniform ellipticity is equivalent to
any of the usual definitions.

Now for the third method.
Suppose $M$ is a \DGP\ with \Ga\ cone $\G$. Choose $A_0\in\G$.  Let $\l_1(A),...,\l_m(A)$ denote
the  $A_0$-eigenvalue functions, and let $\tr \l (A) = \l_1(A) + \cdots +\l_m(A)$ denote the trace.

\Theorem{5.16}  {\sl
The function
$$
M_\d (A) \ =\ \prod_{j=1}^m \left(  \l_j(A) +\d \tr(A)\right)
\eqno{(5.12)}
$$
is a \DGP\ with uniformly elliptic branches $F^\d_k$ containing the branches $F_k$ of $\{M=0\}$.
}
\pf
By Proposition 2.30, $M_\d$ is $A_0$-hyperbolic, and up to scale the 
$A_0$-eigenvalue functions are $\L_k(A) = \l_k(A) +\d\tr(A)$.
Each branch $F^\d_k$  of $\{M_\d =0\}$ is $\G_\d$-monotone by (5.1).
It remains to show that $\Int\cp_{\d'}\ss\G_\d$ for some $\d'>0$.
This follows from (5.10) with $F=\overline{\G_\d}$ since 
$\cp-\{0\} \ss \overline\G -\{0\}\ss \Int \overline{\G_\d} = \G_\d$. (Note that $\tr \l (A) >0$ if $A\in \G$.)\qed

 \Remark{5.17} The three methods of generating new DGPs, and hence new subequations
 for which Theorem 5.5 solves the Dirichlet problem, do not commute with each other.
 Thus combining these results creates a  large number of new \DGP s.\medskip
 
 Additional methods for generating subequations and solving the Dirichlet problem
 are investigated next.
 
 \medskip


 \centerline{\headfont Universal Eigenvalue Subequations}
\medskip

The branches $F_k$ of $\{M=0\}$ are universally determined by requiring that:
$$
\l(A) \ \in \ \E_k
$$
where $\E_k \equiv \{ \lu_k\geq0\} = \{\l\in\bbr^m : $ At least $m-k+1$ of the $\l_j$ are $\geq 0\}$.
Thus  the subset $\E_k\ss \bbr^m $ universally determines the $k$th branch of $\{M=0\}$ for any \DG
 polynomial  of degree $m$ on $\Symn$.
For example, $\E$ defined by $\l_1\geq 0, ... , \l_m\geq0$ is the {\bf universal Monge-Amp\`ere subequation},
and it induces an {\sl $M$-Monge-Amp\`ere subequation}  
for each  \DG polynomial  $M$  of degree $m$.

\Def{5.18} A closed symmetric subset $\E$ of $\bbr^m$, with $\emptyset \neq \E \neq \bbr^m$,
 which is positive-orphant monotone,
will be called a {\bf universal eigenvalue subequation}.

If ,in addition, \smallskip

(a) \ \ $\E$ is convex, then $\E$ will be called a {\bf convex (universal) eigenvalue subequation} 
\smallskip

(b)\ \  $\E$ is $(\Or)_\d$-monotone for some $\d>0$, then $\E$ will be called a {\bf 
uniformly elliptic  (universal) eigenvalue subequation}.

\medskip
The structure of universal eigenvalue subequations  is described  in Proposition 5.25 below.
We leave to the reader a similar result for uniformly elliptic  universal eigenvalue subequations.

\Theorem{5.19}  {\sl
A universal eigenvalue subequation $\E$ on $\bbr^m$ universally determines a 
 subequation $F_\E$ on $\bbr^n$ for each \DG polynomial $M$ on $\Symn$ of degree
 $m$ by setting
 $$
 F_\E\ \equiv\ \{A\in \Symn : \l(A) \in \E\} \ =\ \l^{-1}(E).
 $$
Moreover,
$$
\ft_\E \ =\ F_{\wt \E} 
\and
\Fa_{\E} = F_{{{\overrightarrow {E}}}}
\eqno{(5.13)}
$$
Finally,\smallskip

(a) \ \ If $\E$ is convex, then $F_E$ is a convex subequation.

,\smallskip

(b) \ \ If $\E$ is uniformly elliptic, then $F_E$ is a uniformly elliptic subequation.
}
 \pf
 The important monotonicity result Theorem 3.6 implies that $F_E$ is a subequation
 (utilizing Definitions 5.2 and 5.18).  We leave (5.13) to the reader.  Part (a) is Theorem 3.8  while for (b)  
Theorem 3.7 says that $F_E$ is $\overline \G_\d$-monotone.  As noted in the proof of  Theorem 5.16, 
 this implies that $F_E$ is uniformly elliptic.\qed\medskip

Because of the results proved in [HL$_1$] (stated as Theorem 4.15 here) for each of the subequations
$F_\E$ constructed in Theorem 5.18, the Dirichlet problem can be solved.

\Theorem{5.20}  {\sl
Suppose that $\E$ is a universal eigenvalue subequation on $\bbr^m$, and $M$
is a \DG polynomial on $\rn$ of degree $m$.  Then the Dirichlet problem 
for the subequation $F_\E$ induced by $\E$ and $M$ can be
solved uniquely for all continuous boundary data on any
domain with smooth boundary which is both $\Fa_\E$ and $\Fa_{\wt \E}$ strictly
convex.
}
\medskip


 \centerline{\headfont  Examples of  Universal Eigenvalue Subequations}
\medskip

 The first four examples of universal subequations $\E$ have
  already been discussed along with their induced subequations $F_\E$.
  In these examples $\E$ is a branch of an $e$-hyperbolic polynomial $Q$ on
  $\bbr^m$ whose \Ga \ cone $\overline\G$ contains $\Or$. Taking $M$ to be
  any \DGP\ of degree $m$ and composing $Q$ with the eigenvalue map of $M$,
  one sees that Methods 5.10, 5.11 and 5.13 can now be considered special cases of Theorem 5.18.

\Ex{5.21. (The Branches of the Universal Monge-Amp\`ere Equation)}  These are the subsets
$$
\E_k\ =\ \{\l\in \bbr^m : \l_k^{\uparrow} \geq 0\}  \qquad k=1,...,m.
$$
which are the branches of the universal Monge-Amp\`ere equation $Q(\l) = \l_1\cdots \l_m$.

\Ex{5.22. (The Branches of the Universal p$^{\rm th}$ Elementary Symmetric Function)}  
These subsets $\E\ss\bbr^m$ are defined by requiring that at least $k^*=p-k+1$ of the $e$ -eigenvalues
of 
$$
Q(\l)\ =\ {\sum_{|I|=p}}' \l_{i_1} \cdots \l_{i_p} 
$$ 
are $\geq 0$.

\Ex{5.23. (The   Universal Geometrically $p$-Convex Equation and its Branches)}  
Here the subsets $\E$ of $\bbr^m$ are the branches of 
$$
Q(\l)\ =\ {\prod_{|I|=p}}'   \left(\l_{i_1} +\cdots +\l_{i_p} \right)
$$
defined by requiring at least $r$ of the $p$-fold sums $\l_{i_1} +\cdots +\l_{i_p}$
(these are $p$-times the $e$-eigenvalues of $Q$) to be $\geq0$.

\Ex{5.24. (The Universal $\d$-Uniformly Elliptic Subequation)}  
Fix $\d>0$ and let $\tr \l = \l\cdot e = \l_1+\cdots+\l_m$.    Set
$$
Q_\d (\l)\ =\ \prod_{k=1}^m \left ( \l_k + \d \tr \l\right)
$$
The branches of $\{Q_\d=0\}$ defined by requiring at least $k^*= m-k+1$ of the eigenvalues 
$\L_j =  \l_j + \d \tr \l$, $j=1,...,m$ to be $\geq0$, are universal uniformly elliptic subequations.
Note that as $\d>0$ varies, the principal branches  $\left( \Or \right)_\d$ form a 
conically-fundamental  
neighborhood system for $\Or=\bbr_+^m$.

\medskip

 The following provides a classification of universal subequations.

\Prop{5.25}  {\sl
A subset $\E\ss \bbr^m$ is a universal eigenvalue subequation if and only if 
there exists a $\|\cdot\|^\pm$-Lipschitz symmetric function $f$ on the hyperplane
normal to $e$ with 
}
$$
\E\ =\ \{\l+te : \l\cdot e =0 \ \ {\rm and}\ \ t\geq f(\l)\}
$$
\pf
See Example 3.4.

\Ex{5.26.  (Using Functions to Construct Subequations)} 
For simplicity we first consider the universal case.  Suppose 
$$
E_{c_1} \ =\ \{\l : f_1(\l) \geq c_1\}\and
E_{c_2} \ =\ \{\l : f_2(\l) \geq c_2\}.
$$
The description of the set 
$E_{c} \ =\ \{\l :  f_1(\l)  + f_2(\l) \geq c_2\}$ as $\bigcup_{c_1+c_2 \geq c}\left(  E_{c_1}  \cap E_{c_2} \right)$
is awkward.  Similarly, if $f_1\geq 0$ and $f_2\geq 0$, the set $E_c = \{ \l : f_1(\l)f_2(\l) \geq c\}$
 is awkward to describe.  Consequently in describing examples of universal 
 eigenvalue subequations it is useful to state an obvious result.

\Lemma{5.27} {\sl
Suppose that $f(\l)$ is a continuous symmetric function defined on an open symmetric subset
of $\bbr^m$.  Choose $c\in \bbr$, set $E_c = \{\l : f(\l)\geq c\}$, and assume 
$\emptyset\neq E_c\neq \bbr^m$.
\medskip
\centerline{ If $f$ is non-decreasing in each variable, then $\overline{E_c}$ is a universal subequation.}
\medskip
If, in addition, $f$ is concave, then $\overline{E_c}$ is a convex subequation.
}

\Ex{5.28. (Special Lagrangian Type)}
Suppose $\vf: \bbr\to\bbr$ is strictly increasing.  For each $c\in\bbr$
$$
E_c\ =\ \left  \{\l\in\bbr^m : \sum_{k=1}^m \vf(\l_k) \geq c\right\}
$$
 is a universal eigenvalue subequation.   Taking $\vf:\bbr\to (-{\pi\over 2}, {\pi\over 2})$ 
 to be $\vf(t) = \arctan(t)$, we obtain the {\bf universal special Lagrangian subequation}:
 $$
\E_c \ \equiv \ \left \{\l\in \bbr^m  : \  \sum_{k=1}^m \arctan \l_k \geq c\right\}, 
\qquad c\in \left(-\smfrac {m\pi}2,  \smfrac {m\pi}2\right).
$$

The arctangent function is concave on $t>0$, but convex on $t<0$.  
Consequently, in order for $E_c$ to be convex, the condition $\sum_{k=1}^m \arctan\, \l_k \geq c$ 
must imply that each $\l_k\geq 0$. This proves easily that $E_c$ is convex if  $c\geq {(m-1)\pi\over 2}$.
A detailed calculation due to Yuan [Y, Lemma 2.1] shows that 
$$
E_c\ \ {\rm is\ convex}\quad \iff \quad  c\ \geq \  {(m-2)\pi\over 2} 
\eqno{(5.14)}
$$

Note that the function $f$ in the structure Theorem 5.25 is difficult to describe for this set $E_c$.

\Ex{5.29. (Krylov Type)}
First note that by Theorem 2.9 for any $a$-hyperbolic polynomial $p$ on a vector space $V$
$$
-{p_a'(x)\over p(x)} \  = \ -\sum_{k=1}^m {1\over \l_a^k(x)}\ \ {\rm is \ } \G-{\rm monotone\ on \ }\G,
\eqno{(5.15)}
$$
  i.e., strictly increasing in $t$ on lines $x+tb$ (with $x,b\in\G$ and $t>0$).   Equivalently,
  $$
  {p(x)\over p'(x)} \ \ {\rm is \ } \G-{\rm monotone \ on\ }\G.
  $$
 Since each $p^{(k)}(x)$ is positive on $\G$, this implies that 
 $$
 {p(x)\over p^{(k)}(x)} \ =\  {p(x)\over p^{(1)}(x)}  { p^{(1)}(x)\over p^{(2)}(x)} \cdots  { p^{(k-1)}(x)\over p^{(k)}(x)} 
 $$
 is also $\G$-monotone on $\G$.  Equivalently, each function $-p^{(k)}(x)/p(x)$ is $\G$-monotone on $\G$.
Consequently, we have that for all  $c_0,..., c_{m-1}\geq 0$,  the set 
 $$
F\ =\    \left \{A\in \overline{ \G(p)} : p(A) - \sum_{k=0}^{m-1} c_k r p^{(k)}(A) \ \geq\ 0 \right\}
\eqno{(5.16)}
$$
 is a subequation for any Dirichlet-\Ga\ polynomial of degree $m$ on $\Symn$.

 Obviously there are many other ways to use (5.15) to construct suequations.  In Krylov [K], the Dirichlet problem
 was solved for $F$ defined by (5.16) when $p(A) = \det A$ is the determinant.

 Using the following lemma, one can show that $F$ is a convex subequation.
 
 \Lemma{5.30} {\sl
 Suppose that $p$ is $a$-hyperbolic on $V$.  Then on the \Ga\ cone $\G$ the function
 $p(x)/p_a'(x)$ is concave.
 }
 \pf
 Set $f(x) = p(x)/p_a'(x)$ and note that $f$ is homogeneous of degree one, so that
 (2.3)$'$ can be used.  By Theorem 3.6 it suffices to consider the universal case
 $$
 p(x) \ =\ x_1\cdots x_m \qquad{\rm with\ \ }a\ =\ (1,...,1).
 $$
 Then $f(x) = p(x)/p_a'(x) = 1/(\sum{1\over x_k})$, and 
 $f_x'(y) = ( \sum {{x_k\over y_k}} )/ (\sum {1\over y_k})^2$
 so that (2.3)$'$ becomes 
 $$
 \left(\sum_k {1\over y_k}\right)^2 \ \leq \  \left(\sum_k {1\over x_k}\right) \left(\sum_k {x_k\over y_k^2}\right)
 $$
 which is the Cauchy-Schwartz inequality where $a_k={1\over \sqrt{x_k}}$ and $b_k= {\sqrt{x_k}\over y_k}$.
 \qed
 
 \medskip
 
 The fact that $F$ is a convex subequation is an immediate consequence of the next result.
 \Cor{5.31}
 $$
\log\, {p^{(k)}(x) \over p(x)}, \ \ {\sl and \ therefore \ also\ }\ {p^{(k)}(x) \over p(x)}, \ \ {\sl is\ convex\ on\ }\ \G.
 $$
 \pf Applying the convex increasing function $\vf(t) = -\log\,(-t)$, $t<0$, to the convex function 
 $-{p(x)\over p'(x)}$ shows that  $\log\,{p'(x)\over p(x)}$ is convex.  Hence
 $
 \log\, {p^{(k)}\over p(x)} \ =\ \sum_{j=1}^k  \log\,{ p^{(j)}\over p^{(j-1)}(x)} 
 $
 is convex.  \qed

 \vskip .5in
 

\centerline{\headfont 6.  The Dirichlet Problem on Riemannian Manifolds}
\medskip

Many of the equations discussed above can be carried over to  riemannian manifolds,
and existence and uniqueness results for the Dirichlet problem have been established
in this context in [HL$_2$]. An essential ingredient for obtaining results in this general
setting centers on the concept of a monotonicity cone. The discussion above
fits neatly into this theory since for any $I$-hyperbolic  polynomial on $\Symn$,
the closed \Ga\ cone is a monotonicity cone for each of the branches $F_k =\{\l_k\geq0\}$.
We give here a  brief sketch of these ideas.
 
 Let $X$ be a  riemannian manifold and recall that any $u\in C^2(X)$ has a well-defined
 riemannian Hessian given by 
 $$
( \Hess \, u)(V,W)\ \equiv\ VWu - (\nabla_V W)u
 $$
for vector fields $V$ and $W$,  where $\nabla$ is the Levi-Civita connection.  This defines a
 symmetric 2-form on $TX$, that is, a section $\Hess\, u$ of $\Sym(T^*X)$.  
 
 Suppose now that $F\ss \Sym(T^*X)$ is a closed subset satisfying the positivity condition
 $F+\cp\ss F$ where $\cp \equiv \{A\in \Sym(T^*X): A\geq0\}$.  Then a function $u\in C^2(X)$
 is said to be $F$-subharmonic if $\Hess_x u\in F$ for all $x\in X$.  This notion extends
 to any $u\in \USC(X)$ essentially as in Definition 4.8 above.  With this understood, the  formulation
given in Section 4 of the Dirichlet problem for a domain $\O\ss \ss X$ and the formulation of
strict $F$-convexity for $\bo$ carry over to this setting. 

However, the global nature of the problem requires a further  global hypothesis
before analogues of Theorem 4.15 can be established.  Let $M\ss \Sym(T^*X)$ be a closed
subset such that  $M+\cp\ss M$ and each fibre $M_x$ is a convex cone with vertex at the 
origin.  Then $M$ is called a {\bf monotonicity cone for} $F$ if under point-wise sum,
$$
F+M\ \ss\ F.
$$
 
 The equations of interest here are those which are,in a sense, ``universal'' in geometry.
 We now make this precise.  Fix a closed subset $\bbf \ss \Symn$ with $\bbf+\cp\ss\bbf$ and
suppose that
 $$
 g(\bbf) = \bbf \qquad \fa g\in {\rm O}_n
  \eqno{(6.2)}
 $$
 where O$_n$ acts on $\Symn$ in the standard way.  Then $\bbf$ naturally determines 
 a subequation $F$ on any riemannian manifold $X$ as follows.  Let $e_1,...,e_n$ be 
 an orthonormal tangent frame field defined on an open set $U\ss X$.  This determines
 a trivialization 
 $$
 \Sym(T^*U)\ \harr{\vf_U}{\cong} \ U\times \Symn,
 \eqno{(6.3)}
 $$
 and we define 
 $$
 F\cap  \Sym(T^*U) \ \equiv \ \vf_U^{-1}(U\times \bbf).
  \eqno{(6.4)}
 $$
 By the invariance (6.2)  this set is independent of the 
 choice of orthonormal frame field. The resulting subequation
 $$
 F\ \ss\ \Sym(T^*X) 
 $$
  is said to be {\bf universally determined} by the constant coefficient 
  subequation $\bbf$.
  
  This construction can be generalized.   For $\bbf\ss \Symn$ as above, let
  $$
  G\ =\ G(\bbf) \ \equiv\ \{g\in {\rm O}_n : g(\bbf)=\bbf\},
  $$
  and suppose $X$ is provided with a {\bf  topological $G$ structure}.
  This means that there is given a distinguished covering of $X$ by open sets
  $\{U_\a\}$ with an orthonormal frame field $e^\a = (e^\a_1,...,e^\a_n)$
  on each $U_\a$ so that the each change of framing has values in $G$, i.e., 
  $$
  g_{\a\b} : U_\a\cap U_\b \ \to \ G \ss\ {\rm O}_n.
  $$
  Then in the manner above, {\sl $\bbf$ universally determines a subequation
  $F\ss\J^2(X)$ on any riemannian manifold with  a topological $G$-structure.}

  For example any almost complex manifold with a compatible riemannian
  metric (i.e., for which $J$ is point-wise orthogonal) carries a topological U$_n$-structure.
  
  The following is one of the results proved in [HL$_2$].

  \Theorem{6.1} {\sl
  Let $X$ be a riemannian manifold with a topological $G$-structure.
  Let $\bbf \ss \Symn$ be a subequation with monotonicity cone $\bbm \ss \Symn$
  and suppose both are $G$-invariant. Let $F,M\ss J^2(X)$ be the subequations
  on $X$ universally determined by $\bbf$ and $\bbm$.
 
  Suppose $X$ carries a  $C^2$  strictly $M$-subharmonic function.

Then for each   bounded domain $\O\ss\ss X$ with smooth boundary $\bo$
which is  both $\Fa$ and $\overrightarrow{\ft}$ strictly convex,  
{\rm the Dirichlet problem is uniquely solvable for all continuous boundary data}.
That is, for each $\vf\in C(\bo)$, these exists a unique $u\in C(\ob)$ satisfying:
\medskip

(1)\ \ $u$ is $F$-harmonic on $\O$, and
\medskip
(2)\ \ $u=\vf$ on $\bo$.
}
\medskip

The constructions above give many interesting examples of such universal equations
in geometry.  For any $G$-invariant \DG polynomial, the closed \Ga\ cone $\bbf_1$
is a monotonicity cone for every branch $\bbf_k$ and all branches are $G$-invariant.
All further constructions can be carried over as well.

\vskip .5in


\centerline{\headfont Appendix A.  An Algebraic  Description of the Branches}
\bigskip

Suppose $p(x)$ is $a$-hyperbolic of degree $m$ as in Section 2.
 Each of the branches
$$
F_m\ \supset \ \cdots \ \supset\  F_1
$$ 
 of $\{p=0\}$ can be described by polynomial inequalities involving the elementary 
 symmetric functions.  Normalize so that $p(a)=1$ and recall from (2.20) that
 $$
p(ta+ x) \ =\ \prod_{k=1}^m (t+\l_k(x))\ =\ t^m + \s_1(x) t^{m-1} + \cdots + \s_m(x)
 \eqno{(A.1)}
 $$
 where 
 $$
 \s_k(x) \ =\ \sum_{i_1<\cdots  < i_k} \l_{i_1}(x) \cdots \l_{i_k}(x) 
 $$
 is the $k$th elementary symmetric function of the $a$-eigenvalues of $x$.
 Set $\s_0(x)\equiv 1$.
 Each $\s_k(x)$ is a homogeneous polynomial in $x$ of degree $k$, which is $a$-hyperbolic
 as noted in (3) of Section 2.
 
 Now suppose that 
 $$
 f(t) \ =\ \prod_{k=1}^m (t+\l_k) \ =\ t^m+ \s_1 t^{m-1} + \cdots + \s_m
 $$
 is any monic polynomial with real roots.
Let $E^+$ denote the number of $\l_k\geq0$.
Let Var$(\s)$ denote the total number of strict consecutive sign changes in the signs of 
$\s = (\s_0,\s_1,...,\s_m)$ where $\s_0=1$.  By definition {\sl strict} means first drop all zeros in $\s$ 
and then compute the number of consecutive sign changes in the remaining tuple of real numbers. 
 The classical {\sl Descartes rule of signs} is the statement that
 $$
 E^+ + {\rm Var}(\s)\ =\ m
 \eqno{(A.2)}
 $$

 \Lemma {A.1} {\sl  The $k$th branch $F_k$ of $\{p=0\}$ is given by  $F_k = \{ x :  {\rm Var}(\s(x))\ \leq\ k-1\}$.
 }
\medskip 

 For example, $\bbf_1(M)$ is defined by
 $$
 \s_1(x)\ \geq\ 0,\ \cdots  , \  \s_m(x)\ \geq\ 0
 \eqno{(A.3)}
 $$
 since Var$(\s) \leq 0$ if and only if $\s_1\geq 0$,...,$\s_m\geq 0$.

\pf Using the ordered $a$-eigenvalues of $x$, the $k$th branch  $F_k$ of $\{p=0\}$ is defined by 
$ \l_k(x)\geq0$.
or equivalently by the condition
$$
E^+(x) \ \geq\ m-k+1.
$$ 
By Descarte's rule of signs (A.2), this  is equivalent to Var$(\s(x)) \leq k-1$. \qed

 \Cor{A.2. (Algebraic description of the branches)}  {\sl
 Set $\e = (\e_0,...,\e_m)$ with $\e_0=1$ and each $\e_j=\pm 1$.
 Define $F^\e$ by the conditions
 $$
\e_1 \s_1(x)\ \geq\ 0,\  \cdots \ ,  \ \e_m\s_m(x)\ \geq\ 0.
 \eqno{(A.4)}
 $$
 Then}
 $$
F_k \ =\ \bigcup_{{\rm Var}(\e)\leq k-1}  F^\e.
  \eqno{(A.5)}
 $$
 
 \Ex{A.3}  Consider the polynomial $p(x,y,z) = {1\over 3}(xy+xz+yz)$ in $\bbr^3$ and set 
 ${\bf a}=(1,1,1)$. Then
 $$
 p(t{\bf a}+{\bf x})\ =\ t^2 + \smfrac 2 3 (x+y+z) t  + \smfrac 1 3 (xy+xz+yz)\ =\ t^2+\s_1 t +\s_2.
 $$
One easily computes that  $\{\s_2 \geq 0\}  = C\cup (-C)$ where $C$ is the  convex circular cone
with base $B=\{x+y+z=1\} \cap \{x^2+y^2+z^2\leq 1\}$. 
The hyperplane $\{\s_1=0\}$ divides the set $\{ \s_2 \geq  0\}$ into the two pieces $C$ and $-C$.
The decompositions in Corollary A.2 are
$$
F_1 \ =\ C\ =\ F^{(1,1,1)} \and F_2\ =\ \sim(-C)\ =\ F^{(1,1,1)} \cup F^{(1,-1,-1)}\cup F^{(1,1,-1)}
$$

\vskip .3in


\centerline{\headfont Appendix B.  Gurvits' Inequality}
\bigskip

\Ga's classical inequality and the recent improvement by  Gurvits 
are derived in this appendix. 
We present the proofs in  a way that emphasizes the parallels.
Suppose that  $p$ is an $a$-hyperbolic polynomial of degree $m$,
and let $\l_1(x),...,\l_m(x)$ denote the $a$-eigenvalues of $x\in V$.
  Elementary Property (3) says that
$$
p(a+tb)\ =\ p(a)\prod_{k=1}^m\left ( 1+t\l_a^k(b) \right) \quad \forall \, t\in\bbr, b\in V.
\eqno{(B.1)}
$$
and
$$
p_b'(a)\ =\ p(a)\sum_{k=1}^m \l_a^k(b).
\eqno{(B.2)}
$$

If $\l_1,..., \l_m >0$, then the classical inequality between 
the geometric and arithmetic means says that
$$
\qquad\qquad
\prod_{k=1}^m\left (1+ t \l_k\right ) \ \leq \ \left( 1+ {\sum \l_k\over m} t    \right)^m
\qquad {\rm with\ equality }\iff  \l_1=\cdots =\l_m.
\eqno{(B.3)}
$$
Applying  (B.3)  to (B.1) and using (B.2) to substitute for $\sum \l_k(b)$, proves the following inequality.

\Lemma{B.1} {\sl
Suppose $p$ is an $a$-hyperbolic polynomial of degree $m$
and  $b\in\G(p)$.  Then for $0<t<\infty$
$$
p(a+tb) \ \leq\ p(a)\left(   1 +  {p_b'(a) \over  m \,p(a)} t  \right)^m  
\eqno{(B.4)}
$$
Equality holds if and only if one of the following equivalencies holds:\medskip

(i)\ \ $\l_1(b) \ =\ \cdots \ =\ \l_m(b)$.

\medskip

(ii)\ \  $a$ and $b$ are proportional module the edge $E(p)$.

\medskip

(iii)\ \  $p(sa+tb)\ =\ (\a s+\b t)^m$ for some $\a>0, \b>0$.
 }
\medskip
\noindent
{\bf Proof that (i) $\iff$ (iii).}  Use formula (2.1)$'$ for $p(sa+tb)$.
\medskip
\noindent
{\bf Proof that (i) $\iff$ (ii).}  Recall (Theorem 2.3(a)) that the edge equals the null set of $p$.
By Property (2), $\l_k(b-\mu a) = \l_k(b) -\mu$.  This proves that
$$
\l_k(b) \ = \ \mu \fa k\quad \iff\quad b-\mu a\in\Edge(p).\qquad\qquad \mathqed
\eqno{(B.5)}
$$

\Ga's inequality follows from (B.4) by dividing both sides by $t^m$ and taking infimums,
which occur at $t=\infty$ for both sides.

\Prop{B.2. (\Ga)} {\sl
For all $a,b\in\G(p)$ one has
$$ 
 p(b) \ =\ \inf_{t>0} {1\over t^m} p(a+tb)\ \leq \ p(a)^{1-m}\left ( {1\over m} p_b'(a)\right)^m
\eqno{(B.6)}
$$
with equality as in Lemma B.1.
}
\medskip

Gurvits' inequality follows from (B.4) by dividing both sides by $t$ and taking infimums.
If equality occurs, then the infimums must occur at the same point $t$,
and at this point $t$, equality must also occur in (B.4).  This establishes the equality statement
in the next result.

\Prop{B.3. (Gurvits)} 
{\sl Assume $m\geq 2$. For all $a,b\in \G(p)$
$$
{(m-1)^{m-1} \over m^m} \,\inf_{t>0} {p(a+tb)  \over t}
\ \leq 
 {1\over m}p_b'(a)
\eqno{(B.7)}
$$
with equality the same as in Lemma B.1.}
\pf
By (B.4) for each $t>0$
$$
{p(a+tb)\over t}
\ \leq\ 
p(a) {\left( 1 + {\a\over m}  \right)^m  \over t}  \quad {\rm with} \ \ \a={p_b'(a)\over p(a)}.
\eqno{(B.8)}
$$
The function $f(t)\equiv {1\over t}\left( 1 + {\a\over m}  \right)^m$ on the interval $0<t<\infty$
blows up at the endpoints and has a single critical point at $t={m\over m-1}{1\over \a} = 
{m\over m-1}{ p(a)\over p_b'(a)}$ with critical value
$( {m\over m-1})^m \a = 
( {m\over m-1})^m {p_b'(a)\over p(a)} = 
{m^m \over (m-1)^{m-1} } {1\over p(a)} {1\over m} p_b'(a)$.
This proves  the inequality (B.7). \qed

\Ga's full inequality follows from the special case (B.6) by induction.

\Theorem{ B.4. (\Ga)}  {\sl  Suppose that $p$ is  $a$-hyperbolic   of degree $m$
and $b_1,...,b_m\in \G(p)$.  Then 
$$
p(b_1)^{1\over m} \cdots p(b_m)^{1\over m}\ \leq \   {1\over m!} \, p^{(m)}_{b_1,...,b_m}.
\eqno{(B.9)}
$$
with equality if and only if one of the following equivalencies holds:
\medskip

(i)\ \ $b_1,...,b_m$ are pairwise proportional modulo the edge $E(p)$.

\medskip
 {(ii)}\ \ Edge$(p) \cap W$ has codimension one in $W\equiv \span\{b_1,...,b_m\}$.

\medskip
{(iii)}\  \ $p(t_1b_1+\cdots t_mb_m) \ =\ (\a_1t_1+\cdots +\a_mt_m)^m$ for some $\a_1,...,\a_m>0$.
}

\pf
Induction can be employed because $p_b'$ is also $a$-hyperbolic with $\G(p_b') \supset \G(p)$
and  Edge$(p_b') = $ Edge$(p)$. (See Propositions 2.21 and 2.32.)  
Employing (B.6)$'$ with $a=b_j$, $ 1\leq j\leq m-1$ and with  $b=b_m$ we have
$$
p(b_j)^{m-1\over m} p(b_m)^{1\over m}\ \leq {1\over m} p_{b_m}'(b_j) \qquad {\rm for}\ \ j=1,...,m-1.
$$
Taking the product of both sides over $j=1,...,m-1$ yields
$$
p(b_1)^{m-1\over m} \cdots p(b_m)^{m-1\over m} \ \leq \ {1\over m^{m-1}}p_{b_m}'(b_1) \cdots p_{b_m}'(b_{m-1})
\eqno{(B.10)}
$$
The induction hypothesis applied to the polynomial
$p_{b_m}'$ says that 
$$
p_{b_m}'(b_1)^{1\over m-1} \cdots p_{b_m}'(b_{m-1})^{1\over m-1}\ \leq \ {1 \over (m-1)!}
 p_{b_1,...,b_m}^{(m)}
\eqno{(B.11)}
$$
Taking the $(m-1)^{\rm st}$ root of both sides of (B.10)  and combining it with  (B.11) proves (B.9).  The equality statement follows easily from Lemma B.1.\qed

\Remark{B.5} If $a,b_1,...,b_m \in V$, then the $m$-fold directional derivative in the 
directions $b_1,...,b_m$ is a constant independent of $a$.  Let $\overline p$ denote the completely
polarized form of $p$.  Then
$$
{1\over m!}p^{(m)}_{b_1,...,b_m}
\ =\ 
{1\over m!}{\partial^m  \over \partial t_1 \cdots \partial t_m} p(a+t_1b_1+\cdots+t_mb_m)\biggr|_{t=0}
\ =\ 
\overline p(b_1,...,b_m)
\eqno{(B.12)}
$$

If $q$ is a homogeneous polynomial of degree $k$, then $q_a'(a) = kq(a)$ for any direction $a\in V$.
Setting $b_{k+1} = \cdots = b_m = a\in \G(p)$
in \Ga's inequality (B.9) yields
$$
p(b_1)^{1\over m} \cdots p(b_k)^{1\over m} p(a)^{m-k\over m}\ \leq \   {1\over k!} \, p^{(k)}_{b_1,...,b_k}(a).
\eqno{(B.13)}
$$

Gurvits' full inequality follow from the special case (B.7) by induction.

\Theorem{B.6. (Gurvits)} {\sl
Suppose that $p$ is $a$-hyperbolic and $b_1,...,b_m\in\G(p)$.  Then
$$
{1\over m^m} {\rm Cap}(p)\ \equdef\ 
{1\over m^m}  \inf_{t_1,...,t_m>0} {p\left (t_1b_1+\cdots + t_m b_m\right)  \over t_1\cdots t_m}
\ \leq\ 
{1\over m!} p^{(m)}_{b_1,...,b_m}
\eqno{(B.14)}
$$
with equality  holding exactly as in \Ga's Theorem.
}
\pf
We use Lemma B.3  repeatedly.
$$
{1\over m^m}  \inf_{t_m>0} {p\left (t_1b_1+\cdots + t_m b_m\right)  \over t_1\cdots  t_m}
\ \leq\ 
{1\over m(m-1)^{m-1}}  {p_{b_m}^{(1)}(t_1b_1+\cdots + t_{m-1}b_{m-1} ) \over t_1\cdots t_{m-1}}
$$
and 
$$
\eqalign
{
{1\over m(m-1)^{m-1}}  \inf_{t_{m-1}>0} & {p^{(1)}_{b_m}\left (t_1b_1+\cdots + t_{m-1} b_{m-1}\right)  \over t_1\cdots t_{m-1}} \cr
&\leq\ 
{1\over (m-1)(m-2)^{m-2}}  {p_{b_{m-1}, b_m}^{(2)}(t_1b_1+\cdots +t_{m-1}b_{m-2} ) \over t_1\cdots t_{m-2}}
}
$$
etc., proves that 
$$
{1\over m^m}  \inf_{t_2,...,t_m>0} {p\left (t_1b_1+\cdots + t_m b_m\right)  \over t_1\cdots  t_m}
\ \leq\ 
{1\over m!}  {p_{b_2,...,b_m}^{(m-1)}(t_1b_1 ) \over t_1}
$$
which equals ${1\over m!} p_{b_1,...,b_m}^{(m)}$ since $p_{b_2,...,b_m}^{(m-1)}$ is homogeneous 
of degree 1.\qed

\Remark{B.7. (Capacity)}  The {\bf capacity} of $p$ with respect to 
 $b_1,...,b_m \in   \G(p)$ is defined by 
 $$  
{\rm Cap}_{b_1,...,b_m}(p) \ =\ \inf_{t_1,...,t_m>0} {p(t_1b_1+\cdots + t_mb_m) \over t_1\cdots t_m}
$$
The inductive step above can be written as
$$
{1\over k^k} {\rm Cap}_{b_1,...,b_k} \left ( p^{(m-k)}_{b_{k+1},...,b_m}  \right)
\ \leq\ 
{1\over k(k-1)^{k-1}} {\rm Cap}_{b_1,...,b_{k-1}} \left ( p^{(m-k+1)}_{b_{k},...,b_m}  \right)
\eqno{(B.15)}
$$

\Remark{B.8 (Gurvits Refines/Improves \Ga)}
Suppose $p$ is $a$-hyperbolic and $b_1,...,b_m\in \overline {\G(p)}$.
$$
p(b_1)^{1\over m} \cdots  p(b_m)^{1\over m} \ \leq\  {1\over m^m}\inf_{x_1,...,x_m>0} 
{p(x_1 b_1+\cdots+ x_mb_m) \over x_1\cdots x_m} \ \equdef \  {1\over m^m}  {\rm Cap}(p)
\eqno{(B.16)}
$$
\pf
Using homogeneity and replacing $b_j$ by $x_jb_j$ the inequality becomes
$$
p(b_1)^{1\over m} \cdots  p(b_m)^{1\over m} \ \leq\  {1\over m^m}p(b_1+\cdots+b_m).
\eqno{(B.16)'}
$$
We may assume that $p(b_j) >0$, $j=1,...,m$, i.e., $b_1,...,b_m \in \G(p)$.

Now (B.16)$'$ follows from the classical geometric-arithmetic mean inequality:
$$
p(b_1)^{1\over m} \cdots  p(b_m)^{1\over m} \ \leq\ \left[ {  p(b_1)^{1\over m}+ \cdots + p(b_m)^{1\over m}  \over m}    \right]^m
$$
and the concavity of $p^{1\over m}$:
$$
 p(b_1)^{1\over m}+ \cdots + p(b_m)^{1\over m}  \ \leq \ p(b_1+\cdots+b_m)^{1\over m}
 $$.

\Remark{B.9} Note that the L.H.S. in the \Ga \ inequality vanishes 
if any one of $b_1,...,b_m$ are on the boundary of $\G(p)$.
However, the $(b_1,...,b_m)$-capacity of $p$ occuring in the Gurvits inequality
may not be zero.  See [Gur$_1$] for a more general upper bound result
for the $(b_1,...,b_m)$-capacity which allows for degeneracies.

\vfill\eject

\centerline{\bf References}

\vskip .2in

\noindent
\item{[BGLS]}   H. Bauschke, O. Guler,  A. Lewis and H. Sendov, {\sl   Hyperbolic polynomials and convex analysis},
Canad. J. Math. {\bf 53} no. 3, (2001), 470-488. 

 \smallskip

\noindent
 \item{[CNS]}   L. Caffarelli, L. Nirenberg and J. Spruck,  {\sl
The Dirichlet problem for nonlinear second order elliptic equations, III: 
Functions of the eigenvalues of the Hessian},  Acta Math.
  {\bf 155} (1985),   261-301.

 \smallskip

\noindent
\item{[C]}   M. G. Crandall,  {\sl  Viscosity solutions: a primer},  
pp. 1-43 in ``Viscosity Solutions and Applications''  Ed.'s Dolcetta and Lions, 
SLNM {\bf 1660}, Springer Press, New York, 1997.

 \smallskip

\noindent
\item{[CIL]}   M. G. Crandall, H. Ishii and P. L. Lions {\sl
User's guide to viscosity solutions of second order partial differential equations},  
Bull. Amer. Math. Soc. (N. S.) {\bf 27} (1992), 1-67.

 \smallskip

\noindent
\item{[G$_1$]}   L. G\aa rding, {\sl  Linear hyperbolic differential equations with constant coefficients},
Acta. Math. {\bf 85}   (1951),  2-62.

 \smallskip

\noindent
\item{[G$_2$]}   L. G\aa rding, {\sl  An inequality for hyperbolic polynomials},
 J.  Math.  Mech. {\bf 8}   no. 2 (1959),   957-965.

 \smallskip

\noindent
\item{[Gu]}   Guler, {\sl   Hyperbolic polynomials and interior point methods for convex programming},
 Math.  Oper. Res.  {\bf 22} no. 2 (1997), 350-377.

 \smallskip

\noindent
\item{[Gur$_1$]}   L. Gurvits, {\sl  Hyperbolic polynomials approach to 
Van der Waerden/Schrijver-Valiant like conjectures: sharper bounds, simpler proofs and algorithmic 
applications },
{\sl Proc. 38 ACM Symp. on Theory of Computing (StOC-2006)}, 417-426, ACM, New York, 2006.

 \smallskip

\noindent
\item{[Gur$_2$]}   L. Gurvits, {\sl  Van der Waerden/Schrijver-Valiant like conjectures and stable (aka  hyperbolic) homogeneous polynomials: one theorem for all},
 ArXiv:0711.3496,  May, 2008.

 \smallskip

\item {[HL$_{1}$]}  F. R. Harvey and H. B. Lawson, Jr., {\sl  Dirichlet duality and the non-linear Dirichlet problem},    Comm. on Pure and Applied Math. {\bf 62} (2009), 396-443.

\smallskip

\item {[HL$_{2}$]}  F. R. Harvey and H. B. Lawson, Jr., {\sl  Dirichlet duality and the non-linear Dirichlet problem
on riemannian manifolds},   ArXiv: 0907.1981.

\smallskip

\item {[HL$_{3}$]} F. R. Harvey and H. B. Lawson, Jr., {\sl  Lagrangian plurisubharmonicity and convexity},  Stony Brook Preprint (2010).

\smallskip

\item {[H]}   L. H\"ormander,  {\sl  Lagrangian plurisubharmonicity and convexity},  Stony Brook Preprint (2010).

\smallskip

   \noindent
\item{[K]}    N. V. Krylov,    {\sl  Notions of Convexity},   Progress in Mathematics {\bf 127},
Birkh\'auser, Boston, 1994.
\smallskip

   \noindent
\item{[L]}  P. Lax,    {\sl  Differential equations, difference equations and matrix theory},  
Comm. on Pure and Applied Math. {\bf 11} (1958), 175-194.
\smallskip

   \noindent
\item{[LPR]}  A. S. Lewis, P. A. Parrilo, M. V. Ramana,    {\sl  The Lax conjecture is true},   Proc. Amer. Math. Soc. {\bf 133} (2005), no. 9, 2495--2499.
\smallskip

\noindent
\item{[R]}   J.  Renegar, {\sl Hyperbolic programs, and their derivative relaxations},
Found.  Comput.  Math. {\bf 6} (2006) no. 1, 59-79.
  
 \smallskip

\item {[S]}  Z. Slodkowski, {\sl  The Bremermann-Dirichlet problem for $q$-plurisubharmonic functions},
Ann. Scuola Norm. Sup. Pisa Cl. Sci. (4)  {\bf 11}    (1984),  303-326.

\smallskip

\item {[V]} V. Vinnikov,  {\sl  Selfadjoint determinantal representations of real plane curves},
Math. Ann.   {\bf 296}    (1993),  453-479.

\smallskip

\item {[Y]} Y. Yuan,  {\sl  Global solutions to special lagrangian equations},
Proc. A. M. S.   {\bf 134} no. 5    (2005),  1355-1358.

\smallskip

\end